\documentclass[a4paper, 12pt]{article}
\usepackage{enumerate, theorem}
\usepackage{amsmath, amsfonts, amssymb}
\usepackage[height=22.5cm, width=15cm]{geometry}
\usepackage[all]{xy}
\usepackage{mathrsfs, yfonts}

\newcommand{\parag}[1]{\paragraph{\sc{#1.}} }

\newtheorem{thm}{Th\'eor\`{e}me}[subsection]
\newtheorem{defn}[thm]{D\'efinition}
\newtheorem{cor}[thm]{Corollaire}
\newtheorem{prop}[thm]{Proposition}
\newtheorem{lemma}[thm]{Lemme}

\begin{document}

\title{Un th\'eor\`eme \`a la "Thom-Sebastiani" pour les int\'egrale-fibres.}

\date{28/09/08 ; r\'evis\'ee le 24/09/09.}

\author{Daniel Barlet\footnote{Barlet Daniel, Institut Elie Cartan UMR 7502  \newline
Nancy-Universit\'e, CNRS, INRIA  et  Institut Universitaire de France, \newline
BP 239 - F - 54506 Vandoeuvre-l\`es-Nancy Cedex.France.\newline
e-mail : barlet@iecn.u-nancy.fr}.}

\maketitle

\markright{Thom-Sebastiani...}

\section*{Abstract}

The aim of this article is to prove a Thom-Sebastiani theorem for the asymptotics of the fiber-integrals. This means that we describe the asymptotics of the fiber-integrals of the function \ $f \oplus g : (x,y) \to f(x) + g(y)$ \ on \ $(\mathbb{C}^p\times \mathbb{C}^q, (0,0))$ \ in term of the asymptotics of the fiber-integrals of the holomorphic germs \ $f : (\mathbb{C}^p,0) \to (\mathbb{C},0)$ \ and \ $g :  (\mathbb{C}^q,0) \to (\mathbb{C},0)$. This reduces to compute the asymptotics of a convolution \ $\Phi*\Psi$ \ from the asymptotics of \ $\Phi$ \ and \ $\Psi$ \ modulo smooth terms. \\
To obtain a precise result, giving the non vanishing  of expected singular terms in the asymptotic expansions of the fiber-integrals associated to \ $f\oplus g$, we have to compute the constants coming from the convolution process. We show that they are given by rational fractions of Gamma factors. This enable us to show that these constants do not vanish.\\

\bigskip

\section*{R\'esum\'e.} L'objet de cet article est de d\'emontr\'e un th\'eor\`eme "\`a la Thom-Sebastiani" pour les d\'eveloppements asymptotiques des int\'egrale-fibres des fonctions du type \ $f \oplus g : (x,y) \to f(x) + g(y)$ \ sur \ $(\mathbb{C}^p\times \mathbb{C}^q, (0,0))$ \ en terme des d\'eveloppements asymptotiques des int\'egrale-fibres associ\'ees aux germes holomorphes \ $f : (\mathbb{C}^p,0) \to (\mathbb{C},0)$ \ et \ $g :  (\mathbb{C}^q,0) \to (\mathbb{C},0)$. Ceci se ram\`ene \`a calculer les d\'eveloppements asymptotiques d'une convolution \ $\Phi*\Psi$ \ \`a partir des d\'eveloppements asymptotiques de \ $\Phi$ \ et \ $\Psi$ \ modulo les termes non singuliers.\\
Pour obtenir un r\'esultat pr\'ecis donnant la non nullit\'e des termes singuliers attendus dans les d\'eveloppements asymptotiques des int\'egrale-fibres associ\'ees \`a  \ $f\oplus g$, nous devons calculer les constantes qui apparaissent dans la convolution. Nous montrons qu'elles sont donn\'ees par des fractions rationnelles de facteurs Gamma, ce qui nous permet de montrer qu'elles sont non nulles.\\

AMS Classification (2000) : 32-S-25, 32-S-40, 32-S-50.

\bigskip

Key words : Asymptotic expansions, fiber-integrals, Thom-Sebastiani theorem.\\

\bigskip

Mots clefs : D\'eveloppements asymptotiques, int\'egrale-fibres, th\'eor\`eme de Thom-Sebastiani.

\bigskip

Titre abr\'eg\'e : Thom-Sebastiani ...

\bigskip

Titre anglais : A Thom-Sebastiani theorem for fiber-integrals.

\tableofcontents

\newpage

\section{Introduction.}

Le but du present article est de montrer directement un th\'eor\`eme \`a la "Thom-Sebastiani" (voir [S.-T. 71])  pour les d\'eveloppements asymptotiques des  int\'egrale-fibres d'une fonction de la forme \ $(x,y) \to f(x) + g(y)$ \ o\`u \ $f$ \ et \ $g$ \ sont deux fonctions holomorphes au voisinage de l'origine dans \ $\mathbb{C}^p$ \ et \ $\mathbb{C}^q$ \ respectivement.\\
Il r\'esulte du th\'eor\`eme g\'en\'eral  de d\'eveloppement asymptotique des int\'egrale-fibres (voir [B. 82]), qu'une int\'egrale-fibre, c'est-\`a-dire une fonction de la forme 
$$ s \to F_{\varphi}(s) =  \int_{f=s} \quad \frac{\varphi}{df\wedge d\bar f} $$
o\`u \ $f  : U \to \mathbb{C}$ \ est une fonction holomorphe non constante sur une vari\'et\'e complexe \ $U$ \ de dimension \ $n+1$ \ et \ $\varphi$ \ une \ $(n+1,n+1)-$forme \ $\mathscr{C}^{\infty}$ \ \`a support compact dans \ $U$, peut s'\'ecrire, au voisinage de \ $s = 0$,
$$ F_{\varphi}(s) = \sum_{\alpha \in A, j\in [0,\mu(\alpha)]} \ \theta_{\alpha,j}(s).\vert s\vert^{2\alpha}.(Log\vert s\vert)^j \ + \  \xi(s) $$
o\`u les fonctions \ $\theta_{j,\alpha} $ \ et \ $\xi$ \ sont \ $\mathscr{C}^{\infty}$ \  au voisinage de \ $s = 0$,  o\`u \ $A \subset ]-1,+\infty[ \cap \mathbb{Q}$ \ est un ensemble fini \ et o\`u  \ $\mu : A \to \mathbb{N}$ \ est une application \`a valeurs dans \ $[0,n]$.\\
Nous dirons qu'une fonction  qui admet une telle \'ecriture au voisinage de l'origine admet un {\rm d\'eveloppement asymptotique standard uniforme}\footnote{nous consid\'erons ici des fonctions "uniformes" contrairement aux fonctions holomorphes "multiformes" que l'on obtient dans les int\'egrales "\`a la Malgrange", c'est \`a dire en int\'egrant des formes holomorphes sur des familles horizontales de cycles (voir [M. 74]).} de type \ $(A,\mu)$ \ en \ $s = 0$.\\
En fait, le r\'esultat que nous pr\'esentons consiste essentiellement  \`a montrer la stabilit\'e par convolution de la classe des fonctions admettant \`a l'origine de \ $\mathbb{C}$ \ des "d\'eveloppements asymptotiques  standards uniformes", en pr\'ecisant les exposants et les degr\'es des termes logarithmiques de la convol\'ee \`a partir des informations correspondantes pour les fonctions initiales. \\
Donnons une d\'efinition qui facilitera l'\'enonc\'e de notre r\'esultat.
\begin{defn}
Pour  deux types donn\'es \ $(A,\mu)$ \ et \ $(B,\nu)$ \ nous d\'efinirons le type \ $(A*B, \mu*\nu)$ \ de la fa{\c c}on suivante :
\begin{align*}   
& A*B = \{ \alpha+\beta+1, \alpha \in A, \beta \in B \}.\\
& (\mu*\nu)(\alpha+\beta+1) = \mu(\alpha) + \nu(\beta) \quad {\rm quand} \quad \alpha, \beta \ {\rm et} \ \alpha+\beta+1 \ {\rm ne \ sont\ pas \ dans} \ \mathbb{N}.\\
&  (\mu*\nu)(\alpha+\beta+1) = \mu(\alpha) + \nu(\beta) + 1 \quad {\rm quand} \quad \alpha+\beta+1 \in \mathbb{N}, \alpha \ {\rm et} \ \beta \ {\rm non \ dans} \ \mathbb{N}.\\
& (\mu*\nu)(\alpha+\beta+1) = \mu(\alpha) + \nu(\beta) - 1\quad {\rm quand} \quad \alpha \ {\rm ou} \ \beta \ {\rm (ou \ les \ deux) \ sont \ dans} \  \mathbb{N}.
\end{align*}
\end{defn}

Pr\'ecis\'ement nous montrons le premier th\'eor\`eme suivant :

\begin{thm}[Thom-Sebastiani pour les int\'egrale-fibres.]\label{Th. Seb.}
Soient \ $f$ \ et \ $g$ \ deux fonctions holomorphes au voisinage de l'origine dans \ $\mathbb{C}^{p}$ \ et \ $\mathbb{C}^{q}$ \ respectivement. Supposons que les int\'egrale-fibres de \ $f$ \ et \ $g$ \ admettent des d\'eveloppements asymptotiques  standards uniformes de type respectifs  \ $(A,\mu)$ \ et \ $(B,\nu)$. Alors les int\'egrale-fibres de la fonction \ $f \oplus g : (x,y) \to f(x) + g(y)$ \ admettent des d\'eveloppements asymptotiques  standards uniformes de type \ $(A*B,\mu*\nu)$.
\end{thm}

\bigskip

Mais au del\`a du r\'esultat qualitatif du th\'eor\`eme \ref{Th. Seb.} qui permet de pr\'evoir quel terme peut appara\^itre dans le d\'eveloppement asymptotique standard uniforme de la convolution de deux fonctions admettant de tels d\'eveloppements asymptotiques, il est important de savoir si les termes  "candidats" donn\'es dans ce th\'eor\`eme donnent effectivement des termes non nuls pour la convol\'ee. Ce probl\`eme est assez d\'elicat, et en pr\'esence de plusieurs termes dans les d\'eveloppements asymptotiques des fonctions initiales, des ph\'enom\`enes de compensations peuvent se produire et faire dispara\^itre un terme attendu. Aussi dans la seconde partie de cet article \'etudions-nous en d\'etail le cas de la convolution de deux fonctions admettant un seul terme singulier non nul dans leurs d\'eveloppements asymptotiques respectifs\footnote{en fait on regroupe ensemble les diff\'erentes puissances de \ $Log\vert s\vert$ \ correspondant aux m\^emes exposants pour \ $s$ \ et \ $\bar s$.}. Nous montrons dans ce cas que le terme attendu pour la convolution est effectivement non nul, c'est \`a dire que nous prouvons le th\'eor\`eme pr\'ecis \'enonc\'e ci-dessous.\\
Il n'est pas difficile de se convaincre que pour ce faire on doit calculer pr\'ecisement les constantes qui apparaissent lors de l'op\'eration de convolution, ce qui am\`ene \`a faire toute une s\'erie de calculs qui sont assez fastidieux et pas aussi simples qu'on pourrait le penser \`a priori.\\
Mais, heureusement, les constantes trouv\'ees s'expriment dans tous les cas comme des "fractions rationnelles de facteurs Gammas". On peut donc montrer la non nullit\'e des constantes qui interviennent et donc r\'epondre compl\`etement \`a la question pos\'ee dans le cas consid\'er\'e. \\

Un corollaire combinant la proposition \ref{DA uni} ci-dessous et le th\'eor\`eme pr\'ecis \ref{Precis} montre alors que, quitte \`a bien choisir la forme test et \`a admettre un d\'ecalage des exposants entiers, d\'ecalage born\'e ne d\'ependant que de \ $f,g$ \ et des compacts consid\'er\'es, le terme attendu associ\'e \`a des termes apparaissant effectivement pour \ $f$ \ et \ $g$,  appara\^it effectivement pour la fonction \ $f \oplus g$.\\

On notera que le calcul explicite des constantes permettra \`a l'utilisateur potentiel qui souhaite aller au dela de l'\'enonc\'e de notre "th\'eor\`eme pr\'ecis" et de son corollaire de d\'eterminer dans le cas g\'en\'eral, c'est \`a dire pour un produit de formes test donn\'ees arbitraires, si un ph\'enom\`ene de compensation se produit dans le cas   de la fonction convol\'ee des deux int\'egrale-fibres donn\'ees. Le cas o\`u la somme \ $r + r'$ \ n'admet qu'une seule \'ecriture avec \ $r \in A_f, r' \in B_g$ \  est  probablement simple \`a \'elucider.

\bigskip

Avant de donner l'\'enonc\'e du th\'eor\`eme pr\'ecis, il est important de disposer du r\'esultat suivant qui permet effectivement de satisfaire l'hypoth\`ese de ce th\'eor\`eme.

\bigskip

\begin{prop}\label{DA uni}
Soit \ $\tilde{f} : (\mathbb{C}^{n+1}, 0) \to (\mathbb{C}, 0)$ \ un germe de fonction holomorphe, et soit \ $ f : X \to D$ \ un repr\'esentant de Milnor de \ $\tilde{f}$. Soit \ $K \subset X$ \ un compact ; il existe un entier \ $\kappa$ \ ne d\'ependant que de \ $f$ \ et de \ $K$ \ v\'erifiant la propri\'et\'e suivante:\\
soit \ $\varphi \in \mathscr{C}^{\infty}_K(X,\mathbb{C})^{n,n}$ \ une \ $(n,n)-$forme test telle que la fonction
$$ s \to \int_{f = s} \ \varphi  $$
admette le terme \ $c.\vert s\vert^{2r}.s^m.\bar{s}^{m'}.(Log\vert s\vert )^j $ \ dans son d\'eveloppement asymptotique en \ $s= 0$, avec \ $c \in \mathbb{C}^*$, \ $r \in \mathbb{Q} \cap [0,1[, j \in [0,n]$ \ et \ $(m,m') \in \mathbb{N}^2$, o\`u nous  supposerons l'entier \ $j$ \ maximal pour \ $r,m,m'$ \ donn\'es. Alors pour tout \ $N \geq \kappa + 1 + (m + m')/2 $ \ il existe \ $\varphi_N \in \mathscr{C}^{\infty}_K(X,\mathbb{C})^{n,n}$ \  v\'erifiant
$$  \int_{f = s} \ \varphi_N = \vert s\vert^{2r}.s^{m+\kappa}.\bar{s}^{m'+\kappa}.P_j(Log\vert s\vert ) + \mathcal{O}(\vert s\vert^{2N}) ,$$
o\`u \ $P_j$ \ est un polyn\^ome unitaire de degr\'e \ $j$.
\end{prop}

\bigskip

\begin{thm}[Le th\'eor\`eme de Thom-Sebastiani  pr\'ecis.]\label{Precis}
Soient \ $\tilde{f}$ \ et \ $\tilde{g}$ \ deux germes de fonctions holomorphes au voisinage de l'origine dans \ $\mathbb{C}^{p}$ \ et \ $\mathbb{C}^{q}$ \ respectivement, et notons \ $f : X \to D$ \ et \ $g : Y \to D'$ \ des repr\'esentants de Milnor de ces germes. Soient \ $\varphi \in \mathscr{C}^{\infty}_c(X)^{p,p}$ \ et \ $\psi \in \mathscr{C}^{\infty}_c(Y)^{q,q}$ \ des formes test v\'erifiant les propri\'et\'es suivantes :
\begin{itemize}
\item Il existe des rationnels \ $r,r' \in ]-1,0]$ \ des entiers \ $m,n,m',n',j,j'$ \ et un entier \ $N > m+n+m'+n' +1$ \ tels que l'on ait  
\begin{align*}
& \int_{f=s} \ \frac{\varphi}{df\wedge d\bar f}  = \vert s \vert^{2r}.s^m.\bar s^n.P_j(Log\vert s \vert) + \mathcal{O}(\vert s\vert^{2N}) \\
& \int_{g=s} \ \frac{\varphi}{dg\wedge d\bar g}  = \vert s \vert^{2r'}.s^{m'}.\bar s^{n'}.Q_{j'}(Log\vert s \vert) + \mathcal{O}(\vert s\vert^{2N})
\end{align*}
\end{itemize}
o\`u \ $P_j$ \ et \ $Q_{j'}$ \ sont des polyn\^omes {\bf unitaires} de degr\'es respectifs \ $j$ \ et \ $j'$. On supposera que pour \ $r = 0$ (resp \ $r' = 0$) on a \ $j \geq 1$ (resp. \ $j' \geq 1$).\\
Alors on a, modulo un polyn\^ome \ $\chi_N$ \ de degr\'e (total) \ $\leq 2N-1$ \  l'\'egalit\'e 
\begin{align*}
&  \int_{f(x)+g(y)= s} \frac{\varphi(x)\wedge \psi(y)}{d(f(x)+g(y))\wedge d(\overline{f(x)+g(y)})} = \\
& \\
& \qquad \qquad \qquad \qquad c.\vert s\vert^{2(r+r'+1)}.s^{m+m'}.\bar s^{n+n'}.R(Log\vert s\vert) + \chi_N(s,\bar s) +  \mathcal{O}(\vert s\vert^{2N}) 
\end{align*}
o\`u \ $c$ \ est un nombre complexe {\bf non nul} et o\`u \ $R$ \ est un polyn\^ome {\bf unitaire} dont le degr\'e est d\'etermin\'e de la fa{\c c}on suivante :
\begin{itemize}
\item  si \   $r , r' $ \ et \ $ r +r'+1$ \   sont non nuls le degr\'e de \ $R$ \ est \ $j+j'$.
\item si \ $r$ \ et \ $r'$ \ sont non nuls mais \ $r+r'+1 = 0$ \ le degr\'e de \ $R$ \ est \ $j+j'+1$.
\item si \ $r$ \ ou \ $r'$ \ est nul (ou les deux) le degr\'e de \ $R$ \ est \ $j+j'-1$.
\end{itemize}
\end{thm}

\parag{Remarque} Si on a \ $r = 0 $ \ et \ $j = 0$ \ alors l'int\'egrale fibre de \ $\varphi$ \ ne pr\'esente qu'un terme non singulier (modulo \ $\mathcal{O}(\vert s\vert^{2N})$) \ et donc la convolution des deux int\'egrale-fibres aura  au moins la m\^eme r\'egularit\'e que  l'int\'egrale fibre de \ $\varphi$. On ne peut donc pas esp\'erer de terme singulier dans ce cas. $\hfill \square$

\bigskip

On notera que la valeur pr\'ecise de la constante \ $c$ \ qui d\'epend de \ $r,r',m,m',n,n'$ \ est donn\'ee dans la section 3.

\bigskip

  La d\'emonstration du th\'eor\`eme pr\'ecis est  l'objet de la section 3 et d\'ecoule des corollaires \ref{Cas 1 complet},  \ref{Cas 2 complet} et \ref{Calcul 3}.
  
  \bigskip
  
  Combin\'e avec la proposition \ref{DA uni}, ce th\'eor\`eme donne imm\'ediatement le corollaire suivant qui pr\'ecise "presque compl\`etement" le module des d\'eveloppements asymptotiques de la fonction  \ $ (x,y) \mapsto f(x) + g(y)$ \ \`a partir de ceux de \ $f$ \ et \ $g$.

\bigskip

\begin{cor} Dans la situation du th\'eor\`eme pr\'ec\'edent, fixons des compacts \ $K \subset X$ \ et \ $L \subset Y$ ;  il existe un entier \ $\kappa$, ne d\'ependant que de \ $f,g, K, L$,  tel que, si les int\'egrale-fibres des  formes \ $\varphi \in \mathscr{C}^{\infty}_K(X)^{p,p}$ \ et \ $\psi\in \mathscr{C}^{\infty}_L(Y)^{q,q}$ \ admettent respectivement dans leurs d\'eveloppements asymptotiques les termes
$$ c.\vert s\vert^{2r}.s^m.\bar{s}^{n}.(Log\vert s\vert )^j \quad {\rm resp.} \quad c'.\vert s\vert^{2r'}.s^{m'}.\bar{s}^{n'}.(Log\vert s\vert )^{j'} $$
avec \ $c.c' \not= 0$ \ les entiers \ $j$ \ et \ $j'$ \ \'etant maximaux pour \ $r,m,n$ \ (resp. \ $r', m', n'$) \ donn\'es, alors  il existe une forme test \ $\theta $ \ dans \ $\mathscr{C}^{\infty}_c(X\times Y)^{p+q,p+q}$ \  telle que le d\'eveloppement asymptotique de l'int\'egrale-fibre  
 $$ \int_{f(x)+g(y)= s}\quad  \frac{\theta(x,y)}{d(f(x)+g(y))\wedge d(\overline{f(x)+g(y)})}$$
 contienne le terme
\begin{equation*}
c.c'.\vert s\vert^{2(r+r'+1)}.s^{m+m'+\kappa}.\bar s^{n+n'+\kappa}.R(Log\vert s\vert)
\end{equation*}
o\`u \ $R$ \ est un polyn\^ome {\bf unitaire} dont le degr\'e est d\'etermin\'e de la fa{\c c}on suivante :
\begin{itemize}
\item  si \   $r , r' $ \ et \ $ r +r'+1$ \   sont non nuls le degr\'e de \ $R$ \ est \ $j+j'$.
\item si \ $r$ \ et \ $r'$ \ sont non nuls mais \ $r+r'+1 = 0$ \ le degr\'e de \ $R$ \ est \ $j+j'+1$.
\item si \ $r$ \ ou \ $r'$ \ est nul (ou les deux) le degr\'e de \ $R$ \ est \ $j+j'-1$.
\end{itemize}
\end{cor}

\parag{Remarque} Comme on l'a d\'ej\`a fait remarquer plus haut, des compensations \'eventuelles  entre divers termes des d\'eveloppements asymptotiques ne permettent pas, en g\'en\'eral, d'assurer que la forme test \ $\varphi \wedge \psi$ \ aura dans le d\'eveloppement asymptotique de son int\'egrale-fibre le terme attendu correspondant aux  termes donn\'es dans les d\'eveloppements asymptotiques des int\'egrale-fibres de \ $\varphi$ \ et \ $\psi$.\\
Autrement dit, la forme test \ $\theta$ \ sera, en g\'en\'eral, diff\'erente de \ $\varphi \wedge \psi$, et le d\'ecalge par l'entier \ $\kappa$ \  sera n\'ecessaire.

\bigskip

Une cons\'equence simple de ces r\'esultats est le corollaire suivant,  sur le polyn\^ome de Bernstein  \`a l'origine (voir par exemple [K. 76]) d'une fonction de la forme  \ $(x,y) \mapsto f(x) + g(y)$, qui peut \'egalement  se d\'eduire du th\'eor\`eme de Thom-Sebastiani topologique (voir [Sak.73]) et des r\'esultats de Malgrange-Kashiwara (voir [M.83] ou bien [Bj.93] ch. 6) sur le lien entre la monodromie du complexe des cycles evanescents et les racines du polyn\^ome de Bernstein. \\
Mais, bien s\^ur, les r\'esultats ci-dessus sont bien plus pr\'ecis puisque, non seulement il sont au niveau des termes singuliers des d\'eveloppements asymptotiques, mais aussi puisqu'ils contr\^olent  les d\'ecalages entiers dans les puissances de \ $s$ \ et \ $\bar s$. 

\bigskip

\begin{cor}[Thom-Sebastiani pour le polyn\^ome de Bernstein.]\label{Bern}
Soient \ $f : (\mathbb{C}^p,0) \to \mathbb{C},0)$ \ et \ $g : (\mathbb{C}^q,0) \to (\mathbb{C},0)$ \ deux germes non constants  de fonctions holomorphes. Alors toute racine du polyn\^ome de Bernstein \`a l'origine du germe \ $f \oplus g : (x,y) \to f(x) + g(y)$ \ est modulo \ $\mathbb{Z}$ \ somme d'une racine du polyn\^ome de Bernstein de \ $f$ \ en \ $0$ \ et d'une racine du polyn\^ome de Bernstein de \ $g$ \ en \ $0$.\\
R\'eciproquement, si \ $\alpha$ \ et \ $\beta$ \ sont des racines respectivement des polyn\^omes de Bernstein de \ $f$ \ et \ $g$ \ \`a l'origine, il existe une racine \ $\gamma$ \ du du polyn\^ome de Bernstein de \ $ f\oplus g : (x,y) \to f(x) + g(y)$ \ qui est congrue modulo \ $\mathbb{Z}$ \ \`a \ $\alpha + \beta$.
\end{cor} 

\bigskip 

\section{Le th\'eor\`eme de convolution.}

\subsection{Utilisation du th\'eor\`eme de Fubini.}

Commen{\c c}ons par montrer comment le th\'eor\`eme de Fubini permet de ramener notre probl\`eme \`a la d\'etermination du d\'eveloppement asymptotique \`a l'origine d'une convolution.

   \begin{prop}\label{T-S}
  Soient \ $f : X \to D$ \ et \ $g : Y \to D$ \ des repr\'esentants de Milnor de deux germes de fonctions holomorphes \`a l'origine de \ $\mathbb{C}^p$ \ et \ $\mathbb{C}^q$ \ respectivement. Pour \ $\rho$ \ (resp. \ $\sigma$ ) une fonction \ $\mathscr{C}^{\infty}$ \ \`a support compact dans \ $X$ \ (resp. dans \ $Y$) \ notons \  $\Phi_{f,\rho}$ \ (resp. \ $\Psi_{g,\sigma}$ ) l'int\'egrale fibre
  \begin{align*}
  &  \Phi_{f,\rho}(s) : = \int_{f(x) = s} \ \rho(x).\frac{dx\wedge d\bar x}{df\wedge d\bar f} \quad {\rm ( et \ respectivement) } \\
  & \qquad\\
  &  \Psi_{g,\sigma}(t) : = \int_{g(y) = t} \ \sigma(y).\frac{dy\wedge d\bar y}{dg\wedge d\bar g}
  \end{align*}
  Alors on a l'\'egalit\'e, pour \ $s\in \mathbb{C}$ :
  $$ \int_{f(x) + g(y) = s} \ \rho(x).\sigma(y).\frac{dx\wedge d\bar x\wedge dy\wedge d\bar y}{d(f+g)\wedge d\overline{(f + g)}} = \int_{\mathbb{C}} \ \Phi_{f,\rho}(s-u).\Psi_{g,\sigma}(u).du\wedge d\bar u .$$ 
  \end{prop}
  
  \parag{Preuve} Commen{\c c}ons par remarquer que, d'apr\`es le Nullstellensatz, il existe un entier \ $N$ \ en une $(p-1,p-1)-$forme \ $\omega$ \ de classe \  $ \mathscr{C}^{\infty}$ \ au voisinage de \ $Supp(\rho)$ \ dans \ $\mathbb{C}^p$ \ v\'erifiant 
  $$\omega\wedge df\wedge d\bar f = \vert f(x)\vert^{2N}.dx\wedge d\bar x .$$
  On a alors, au voisinage de \ $Supp(\rho)\times Supp(\sigma)$, l'\'egalit\'e
  $$ \omega\wedge dy\wedge d\bar y\wedge d(f+g)\wedge d\overline{(f+g)} = \vert f(x)\vert^{2N}.dx\wedge d\bar x\wedge dy\wedge d\bar y .$$
  On peut alors \'ecrire, toujours au voisinage de \ $Supp(\rho)\times Supp(\sigma)$, sur la fibre \ $\{f(x) + g(y) = s \}$ :
  $$  \frac{dx\wedge d\bar x\wedge dy\wedge d\bar y}{d(f+g)\wedge d\overline{(f+g)}} = \frac{dx\wedge d\bar x}{df\wedge d\bar f}\wedge dy\wedge d\bar y .$$
  Alors le Th\'eor\`eme de Fubini donne
  \begin{align*}
  & \int_{f(x) + g(y) = s} \ \rho(x).\sigma(y).\frac{dx\wedge d\bar x\wedge dy\wedge d\bar y}{d(f+g)\wedge d\overline{(f + g)}} \quad = \\
  &\quad \\
  &   \int_{y \in  \mathbb{C}^q} \sigma(y).dy\wedge d\bar y.\Big( \int_{f(x)= s-g(y)} \rho(x). \frac{dx\wedge d\bar x}{df\wedge d\bar f}\, \Big) \quad  = \\
  & \quad \\
  &   \int_{u \in \mathbb{C}} \Phi_{f,\rho}(s-u).\Big( \int_{g(y) = u} \sigma(y).\frac{dy\wedge d\bar y}{dg\wedge d\bar g}\, \Big).du\wedge d\bar u \quad  = \\
  & \quad \\
  &   \int_{u \in \mathbb{C}} \Phi_{f,\rho}(s-u).\Psi_{g,\sigma}(u).du\wedge d\bar u \qquad \qquad \qquad \qquad \qquad \blacksquare
  \end{align*}
  
  \bigskip

  \subsection{Le th\'eor\`eme de convolution.}

Montrons maintenant notre th\'eor\`eme de d\'eveloppement asymptotique pour une convolution.

\begin{thm}\label{T-S}
Consid\'erons deux nombres r\'eels \ $\alpha$ \ et \ $\beta$ \ strictement plus grands que \ $-1$ \ et soient \ $\theta$ \ et \ $\eta$ \ deux fonctions dans \ $\mathscr{C}^{\infty}_c(\mathbb{C})$, et posons, pour \ $s \in \mathbb{C}^*$,
$$ F_{\alpha,\beta,j,k}(s) : = \frac{1}{2i\pi}.\int_{\mathbb{C}} \theta(s-u).\eta(u).\vert s-u\vert^{2\alpha}.(Log\vert s-u\vert^2)^j.\vert u\vert^{2\beta}.(Log\vert u\vert^2)^k.du\wedge d\bar u .$$
Alors il existe des fonctions \ $(\zeta_l)_{l \in [0,L]}$ \ et \ $\xi$ \ dans  \ $\mathscr{C}^{\infty}_c(\mathbb{C})$ \ telles que l'on ait
$$ F_{\alpha,\beta,j,k}(s) = \sum_{l=0}^L \ \zeta_l(s).\vert s\vert^{2(\alpha+\beta+1)}.(Log\vert s\vert^2)^l + \xi(s) $$
o\`u l'on a 
\begin{itemize}
\item \ $L = j+k $ \ si \ $\alpha, \beta ,\alpha+\beta+1 \not\in \mathbb{N}$, 
\item \ $L = j+k+1$ \ si \ $\alpha, \beta \not\in \mathbb{N}$ \ et si \ $\alpha+\beta+1 \in \mathbb{N}$,
\item  \ $L = j+k-1$ \ si  \ $\alpha$ \ ou \ $\beta$ \ ou les deux sont entiers. 
\end{itemize}
De plus, on a \ $L = -1$ \ d\`es que \ $\alpha \in \mathbb{N}$ \ et \ $j = 0$ \ ou, sym\'etriquement, d\`es que \ $\beta \in \mathbb{N}$ \ et \ $k = 0$.
\end{thm}

\parag{D\'emonstration} Nous la ferons  en plusieurs \'etapes.
\parag{Premi\`ere \'etape} Remarquons d\'ej\`a que les conditions \ $\alpha > -1$ \ et \ $\beta > -1$ \ assurent que, pour \ $s \not= 0$ \ la fonction int\'egr\'ee est localement int\'egrable. De plus le support compact de la fonction \ $\eta$ \ montre l'int\'egrabilit\'e globale. La fonction est donc bien d\'efinie sur \ $\mathbb{C}^*$. \\
Le changement de variable \ $v = s - u$ \ dans l'int\'egrale d\'efinissant la fonction \ $F_{\alpha,\beta,j,k}$ \ montre que les triplets \ $(\theta, \alpha,j)$ \ et \ $(\eta, \beta,k)$ \ jouent des r\^oles sym\'etriques.\\
Par ailleurs il est clair que pour \ $\alpha \in \mathbb{N}$ \ et \ $j=0$ \ on a une convolution entre une fonction \ $\mathscr{C}^{\infty}_c$ \ est une fonction int\'egrable \`a support compact . On obtient donc une fonction \ $\mathscr{C}^{\infty}_c$, et de m\^eme pour \ $\beta \in \mathbb{N}$ \ et \ $k = 0$.

\parag{Seconde \'etape} Fixons \ $\varepsilon > 0 $. L'int\'egrale
$$ \frac{1}{2i\pi}.\int_{\vert u \vert \geq \varepsilon} \ \theta(s-u).\eta(u).\vert s-u\vert^{2\alpha}.(Log\vert s-u\vert^2)^j.\vert u\vert^{2\beta}.(Log\vert u\vert^2)^k.du\wedge d\bar u$$
donne manifestement une fonction \ $\mathscr{C}^{\infty}$ \ sur l'ouvert \ $ \{ \vert s \vert < \varepsilon/2\} $, puisque maintenant la fonction int\'egr\'ee est  \ $\mathscr{C}^{\infty}$ \ en \ $s$ \ et major\'ee par une fonction int\'egrable fixe ainsi que chacune de ses d\'eriv\'ee partielles en \ $s, \bar s$. \\
Pour prouver l'existence du d\'eveloppement asymptotique \`a l'origine de la fonction \ $F_{\alpha,\beta,j,k}$ \ il suffit donc de le faire pour la fonction obtenue en int\'egrant seulement sur le disque \ $\{ \vert u \vert \leq \varepsilon \}$.

\parag{Troisi\`eme \'etape} Montrons maintenant que si la fonction \ $\theta$ \ est plate \`a l'ordre \ $N \geq 1$ \ \`a l'origine, alors la fonction \ $F_{\alpha,\beta,j,k}$ \ est de classe \ $\mathscr{C}^{N-1}$ \ au voisinage de l'origine. Comme notre hypoth\`ese permet d'\'ecrire
$$ \theta(v) = \sum_{i=0}^{N+1} \ v^i.\bar v^{N-i+1}.\theta_i(v) $$
o\`u les fonctions \ $\theta_i$ \ sont dans \ $\mathscr{C}^{\infty}_c(\mathbb{C})$, il s'agit de voir que, pour chaque \ $i \in [0,N+1]$,  la fonction \ $s \mapsto G_i(s)$ \ d\'efinie par l'int\'egrale
$$    \int_{\vert u\vert \leq \varepsilon} \  (s-u)^i.\overline{(s-u)}^{N-i+1}.\theta_i(s-u).\eta(u).\vert s-u\vert^{2\alpha}.(Log\vert s-u\vert^2)^j.\vert u\vert^{2\beta}.(Log\vert u\vert^2)^k.du\wedge d\bar u $$
est de classe \ $\mathscr{C}^{N-1}$ \ au voisinage de l'origine. Montrons ceci par r\'ecurrence sur \ $N \geq 1$. Pour \ $N = 1$ \ il s'agit de voir que les fonctions \ $G_0, G_1$ \ et \ $G_2$ \ sont continues au voisinage de l'origine. Ceci r\'esulte imm\'ediatement du fait que la fonction
$$ (s-u)^i.\overline{(s-u)}^{2-i}.\theta_j(s-u).\vert s-u\vert^{2\alpha}.(Log\vert s-u\vert^2)^j $$
est continue en \ $(s,u)$ \ pour \ $i = 0,1,2$ \ puisque \ $\alpha > -1$, et de l'int\'egrabilit\'e de la fonction \ $u \mapsto \eta(u).\vert u\vert^{2\beta}.(Log\vert u\vert^2)^k$.\\
Supposons l'assertion montr\'ee pour \ $N-1 \geq 1$ \ et montrons-l\`a pour \ $N$.\\
On a, d'apr\`es le th\'eor\`eme de d\'erivation de Lebesgue, 
\begin{align*}
& \frac{\partial G_i}{\partial s}(s) =  (i+\alpha).\int_{\vert u\vert \leq \varepsilon}  (s-u)^{i-1}.\overline{(s-u)}^{N-i+1}.\theta_i(s-u).\eta(u).R_{\alpha,\beta,j.k}(s,u).du\wedge d\bar u  + \\
&   \int_{\vert u\vert \leq \varepsilon}  \ (s-u)^i.\overline{(s-u)}^{N-i+1}.\frac{\partial \theta_i}{\partial v}(s-u).\eta(u).R_{\alpha,\beta,j,k}(s,u).du\wedge d\bar u + \\
& j.\int_{\vert u\vert \leq \varepsilon}  (s-u)^{i-1}.\overline{(s-u)}^{N-i+1}.\theta_i(s-u).\eta(u).R_{\alpha,\beta,j-1,k}(s,u).du\wedge d\bar u 
\end{align*}
o\`u l'on a pos\'e, pour simplifier l'\'ecriture,
\begin{align*}
& R_{\alpha,\beta,j,k}(s,u) = \vert s-u\vert^{2\alpha}.(Log\vert s-u\vert^2)^j.\vert u\vert^{2\beta}.(Log\vert u\vert^2)^k  
\end{align*}

On  constate alors que, comme les fonctions
$$ v^{i-1}.\bar v^{N-i+1}.\theta_i(v) \quad {\rm et} \quad v^i.\bar v^{N-i+1}.\frac{\partial \theta_i}{\partial v}(v) $$
sont plates \`a l'ordre \ $N-1$ \ \`a l'origine, l'hypoth\`ese de r\'ecurrence donne que  \ $ \frac{\partial G_i}{\partial s}$ \ est de classe \ $\mathscr{C}^{N-2}$ \ au voisinage de l'origine. Comme il en est de m\^eme pour \ $\frac{\partial G_i}{\partial \bar s}$, par un calcul analogue \`a ce que l'on vient de voir, ceci ach\`eve la preuve de notre assertion.

\parag{Quatri\`eme \'etape} Etudions maintenant les fonctions suivantes :
$$ G_{\alpha,\beta,j,k}^{p,q}(s) : = \int_{\vert u\vert \leq 1} u^p.\bar u^q.\vert s-u\vert^{2\alpha}.(Log\vert s-u\vert)^j.\vert u\vert^{2\beta}(Log\vert u\vert)^k.du\wedge d\bar u $$
o\`u \ $j,k,p, q$ \ sont des entiers et o\`u \ $\alpha$ \ et \ $\beta$ \ sont des r\'eels strictement plus grands que \ $-1$. On va \'etablir la proposition suivante :

\begin{prop}
La fonction \ $G_{\alpha,\beta,j,k}^{p,q}$ \  est de la forme suivante:
\begin{itemize}
\item pour \ $p \geq q$ :
$$\sum_{l=0}^{L}  c_l.\vert s\vert^{2(\alpha+\beta+1)}.s^p.\bar s^{q}.(Log\vert s\vert)^l \  +  \ s^{p-q}.\Phi(\vert s\vert^2) $$
\item pour \ $p \leq q $ :
$$ \sum_{l=0}^{L}  c_l.\vert s\vert^{2(\alpha+\beta+1)}.s^p.\bar s^{q}.(Log\vert s\vert)^l \  \ + \  \bar s^{q-p}.\Phi(\vert s\vert^2) $$
\end{itemize}
o\`u \ $\Phi$ \ est une fonction analytique au voisinage de l'origine qui d\'epend de \ $\alpha,\beta,j,k,p,q$, o\`u les constantes \ $c_l$ \ d\'ependent \'egalement de \ $\alpha,\beta,j,k,p,q$, et o\`u l'on a
\begin{itemize}
\item \ $L = j + k$ \ quand \ $\alpha, \beta, \alpha + \beta + 1 \not\in \mathbb{N}$,
\item  \  $L = j + k + 1$ \ quand \ $\alpha + \beta + 1 \in \mathbb{N}$ \ et \ $\alpha,\beta \not\in \mathbb{N}$,
\item \ $ L = j + k - 1$ \ quand \ $\alpha$ \ ou \ $\beta$ \ ou les deux sont entiers. 
\end{itemize}
\end{prop}

\bigskip

On notera que ceci implique que la fonction \ $G_{\alpha,\beta,j,k}^{p,q}$ \ est \ $\mathscr{C}^{\infty}$ \ en dehors de l'origine et admet, quand \ $s \to 0$ \ un d\'eveloppement asymptotique avec des termes "singuliers"\footnote{c'est-\`a-dire qui ne sont pas dans \ $\mathbb{C}[[s,\bar s]]$.} \'eventuels bien pr\'ecis, et que ce d\'eveloppement asymptotique est ind\'efiniment d\'erivable, c'est-\`a-dire que toute d\'eriv\'ee partielle de \ $G_{\alpha,\beta,j,k}^{p,q}$ \ admet quand \ $s \to 0$ le d\'eveloppement asymptotique obtenu en effectuant formellement les m\^emes d\'eriv\'ees partielles sur le d\'eveloppement asymptotique trouv\'e que celles effectu\'ees sur \ $G_{\alpha,\beta,j,k}^{p,q}$.

\parag{Preuve} Posons \ $u = s.t $ \ pour \ $s \not= 0$. Alors \ $G_{\alpha,\beta,j,k}^{p,q}(s)\big/ s^p.\bar s^q.\vert s\vert^{2(\alpha+\beta+1)}$ \ est donn\'e par l'int\'egrale :
$$ \int_{\vert t\vert \leq 1/\vert s\vert} \ t^p.\bar t^q.\vert 1-t\vert^{2\alpha}.\vert t\vert^{2\beta}.(Log\vert s.(1-t)\vert)^j.(Log\vert s.t\vert)^k.dt\wedge d\bar t .$$
Commen{\c c}ons par supposer que \ $\alpha+\beta+1 \not\in \mathbb{N}$.\\
Comme l'int\'egrale pour \ $\vert t\vert \leq 3$ \ est un polyn\^ome en \ $Log\vert s\vert$ \ de degr\'e au plus \ $j+k$ \ dont les coefficients d\'ependent de fa{\c c}on \ $\mathscr{C}^{\infty}$ \ de \ $(\alpha,\beta) \in ]-1,+\infty[^2$, nous pouvons nous contenter d'\'etudier la fonction \ $ \Gamma_{\alpha,\beta,j,k}^{p,q}(s) \big/ s^p.\bar s^q.\vert s\vert^{2(\alpha+\beta+1)}$ \  donn\'ee par l'int\'egrale :
$$  \int_{3 \leq \vert t\vert \leq 1/\vert s\vert} \ t^p.\bar t^q.\vert 1-t\vert^{2\alpha}.\vert t\vert^{2\beta}.(Log\vert s.(1-t)\vert)^j.(Log\vert s.t\vert)^k.dt\wedge d\bar t .$$
Ecrivons 
\begin{align*}
& ( Log\vert s.t\vert)^k = (Log\vert s\vert + Log\vert t\vert)^k \quad\quad  {\rm et}\\
& (Log\vert s.(1-t)\vert)^j = (Log\vert s\vert + Log\vert t\vert + Log\vert 1-1/t\vert)^j 
\end{align*}
et d\'eveloppons notre int\'egrale par la formule du bin\^ome. On obtient ainsi un polyn\^ome de degr\'e \ $\leq j+k$ \ en \ $Log\vert s\vert$ \ et le coefficient de \ $(Log\vert s\vert)^{j+k-h}$ \ est, \`a une combinaison lin\'eaire \`a coefficients constants de  fonctions de la forme
$$\int_{3 \leq \vert t\vert \leq 1/\vert s\vert} \ t^p.\bar t^q.\vert 1-t\vert^{2\alpha}.\vert t\vert^{2\beta}.(Log\vert 1-1/t\vert)^{h_1}.(Log\vert t\vert)^{h_2}.dt\wedge d\bar t $$
o\`u \ $h_1 + h_2 = h$.

\bigskip

En coordonn\'ees polaires \ $t = \rho.e^{i\varphi}$ \ l'int\'egrale pr\'ec\'edente donne, \`a une constante non nulle pr\`es,
$$ \int_3^{1/\vert s\vert} \rho^{p+q+2(\alpha+\beta+1)}.(Log\, \rho)^{h_2}.\frac{d\rho}{\rho}.\int_0^{2\pi} \vert 1- \frac{e^{-i\varphi}}{\rho}\vert^{2\alpha}.(Log\vert 1- \frac{e^{-i\varphi}}{\rho}\vert)^{h_1}.e^{i(p-q)\varphi}.d\varphi .$$

Comme on a \ $\rho \geq 3$ \ on a un d\'eveloppement en s\'erie qui est  normalement convergeant
$$ \vert 1- \frac{e^{-i\varphi}}{\rho}\vert^{2\alpha}.(Log\vert 1- \frac{e^{-i\varphi}}{\rho}\vert)^{h_1} = \sum_{\nu = 0}^{\infty} \gamma_{\alpha,\nu}^{h_1}(cos \varphi).\rho^{-\nu} $$
ce qui donne que \ $\Gamma_{\alpha,\beta,j,k}^{p,q}(s) \big{/}  s^p.\bar s^q.\vert s\vert^{2(\alpha+\beta+1)} $ \ est combinaison lin\'eaire \`a coefficients constants quand \ $h_2 $ \ d\'ecrit \ $[0,j+k]$ \ de termes tels que
\begin{align*}
&  \sum_{\nu = 0}^{\infty} \ C_{\alpha,\nu}^{h_1}. \int_3^{1/\vert s\vert} \rho^{p+q+2(\alpha+\beta+1)-\nu}.(Log \,\rho)^{h_2}.\frac{d\rho}{\rho} \tag{*}
\end{align*}
 o\`u \ $C_{\alpha,\nu}^{h_1} : = \int_0^{2\pi} \ \gamma_{\alpha,\nu}^{h_1}(cos \varphi).e^{i.(p-q).\varphi}.d\varphi $.\\
On notera  que le polyn\^ome en \ $cos(\varphi), \gamma_{\alpha,\nu}^{h_1}(cos(\varphi))$, est une fonction \ $\mathscr{C}^{\infty}$ \ de \\
 $\alpha \in ]-1,+\infty[$, et donc que \ $C_{\alpha,\nu}^{h_1}$ \ d\'epend \'egalement de fa{\c c}on \ $\mathscr{C}^{\infty}$ \ de \ $\alpha \in ]-1,+\infty[$.
 
 \bigskip

Remarquons maintenant que la fonction
$$ H_{\alpha,h_1}^{p,q}(x) = \int_0^{2\pi} \vert 1- x.e^{-i\varphi}\vert^{2\alpha}.(Log\vert 1- x.e^{-i\varphi}\vert)^{h_1} .e^{i(p-q)\varphi}.d\varphi $$
dont le d\'eveloppement en s\'erie \`a l'origine est donn\'e par
$$ H_{\alpha,h_1}^{p,q}(x) = \sum_{\nu = 0}^{\infty} \ C_{\alpha,\nu}^{h_1}.x^{\nu} $$
est \ $(-1)^{p+q}-$paire. Ceci se voit facilement en changeant \ $x$ \ en \ $-x$ \ et \ $\varphi$ \ en \ $\varphi + \pi$.\\ 
Ceci montre que \ $C_{\alpha,\nu}^{h_1}$ \ est non nul seulement quand les entiers \ $p + q$ \ et \ $\nu$ \ ont m\^eme parit\'e.\\
 On en d\'eduit que, sous notre hypoth\`ese \ $\alpha+\beta + 1 \not\in \mathbb{N}$, l'exposant \\ $ 2(\alpha+\beta+1) + p + q  - \nu = 0 $ \ n'appara\^it pas dans le d\'eveloppement \ $(^*)$.
 
 \bigskip
 
 On a donc, pour \ $\alpha + \beta +1 \not\in \mathbb{N}$ :

\begin{equation*}
 \Gamma_{\alpha,\beta,j,k}^{p,q}(s) = s^p.\bar s^q.\vert s\vert^{2(\alpha+\beta+1)}.P_{j,k}(Log\vert s\vert)  + (\frac{s}{\vert s\vert})^p.(\frac{\bar s}{\vert s\vert})^q.\Phi_{j,k}(\vert s\vert)  \tag{**}
 \end{equation*}
o\`u \ $P_{j,k}$ \ est un polyn\^ome de degr\'e \ $\leq j+k$, dont les coefficients d\'ependent de \ $(\alpha,\beta)$ \ de fa{\c c}on \ $\mathscr{C}^{\infty}$ \ ainsi que de \ $p$ \ et \ $q$, et o\`u \ $\Phi_{j,k}$ \ (qui d\'epend \'egalement de \ $\alpha,\beta,p,q$)  est analytique au voisinage de l'origine.\\
On remarquera que le calcul explicite des coefficients du d\'eveloppement en s\'erie de la fonction \ $\Phi_{j,k}$ \ montre qu'elle d\'epend de fa{\c c}on \ $\mathscr{C}^{\infty}$ \ de \ $(\alpha,\beta) \in (]-1, +\infty[)^2$, $\alpha+ \beta+1 \not\in \mathbb{N}$.

\bigskip

Montrons maintenant par r\'ecurrence sur \ $j+k$ \ que la fonction 
$$ s \mapsto (\frac{s}{\vert s\vert})^p.(\frac{\bar s}{\vert s\vert})^q.\Phi_{j,k}(\vert s\vert) $$
est \ $\mathscr{C}^{\infty}$ \ en \ $s$ \ au voisinage de l'origine.\\

Pour \ $\lambda \in \mathbb{C}^*$, en posant \ $u = \lambda.v$ \ dans l'int\'egrale qui d\'efinit la fonction \ $G_{\alpha,\beta,0,0}^{p,q}$, on obtient : 
$$ G_{\alpha,\beta,0,0}^{p,q}(\lambda.s) = \lambda^p.\bar \lambda^q.\vert \lambda\vert^{2(\alpha+\beta+1)}.\int_{\vert v\vert \leq 1/\vert \lambda\vert} v^p.\bar v^q.\vert s-v\vert^{2\alpha}.\vert v\vert^{2\beta}.dv\wedge d\bar v $$
ce qui donne, puisque la fonction
$$ s \mapsto \int_{1 \leq \vert v\vert \leq 1/\vert \lambda\vert} v^p.\bar v^q.\vert s-v\vert^{2\alpha}.\vert v\vert^{2\beta}.dv\wedge d\bar v $$
est \ $\mathscr{C}^{\infty}$ \ au voisinage de \ $s = 0$, 
$$ G_{\alpha,\beta,0,0}^{p,q}(\lambda.s) -   \lambda^p.\bar \lambda^q.\vert \lambda\vert^{2(\alpha+\beta+1)}.G_{\alpha,\beta,0,0}^{p,q}(s) \in \mathscr{C}^{\infty}. $$
On aura donc, pour tout \ $\lambda \in \mathbb{C}^*$ \ fix\'e,
\begin{equation*}
 \frac{1}{\vert \lambda\vert^{p+q}}.\lambda^p.\bar\lambda^q. (\frac{s}{\vert s\vert})^p.(\frac{\bar s}{\vert s\vert})^q.\big{[}\Phi_{0,0}(\vert \lambda.s\vert) - \vert \lambda\vert^{2(\alpha+\beta+1)+p +q}.\Phi_{0,0}(\vert s\vert)\big{]}  \in \mathbb{C}[[s,\bar s]]. \tag{$^{***}$}
 \end{equation*}
\'Ecrivons le d\'eveloppement \`a l'origine de la fonction analytique
$$ \Phi_{0,0}(x) = \sum_{m =0}^{\infty} \ c_m.x^m .$$
La relation \ $(^{***})$ \ donne alors que
$$  \frac{1}{\vert \lambda\vert^{p+q}}.\lambda^p.\bar\lambda^q. (\frac{s}{\vert s\vert})^p.(\frac{\bar s}{\vert s\vert})^q.\Big[\sum_{m=0}^{\infty} \ c_m.(\vert \lambda\vert^m -  \vert \lambda\vert^{2(\alpha+\beta+1)+p +q}).\vert s\vert^m \Big] \in \mathbb{C}[[s,\bar s]].$$
Comme on a suppos\'e  \ $(\alpha +\beta + 1) \not\in \mathbb{N}$, et que \ $\Phi_{0,0}$ \ est \ $(-1)^{p+q}-$paire, on en d\'eduit que pour chaque \ $m \in \mathbb{N}$ \  tel que \ $c_m \not= 0$, on a
$$ (\frac{s}{\vert s\vert})^p.(\frac{\bar s}{\vert s\vert})^q.\vert s\vert^m \in \mathbb{C}[[s,\bar s]]$$
ce qui prouve notre assertion pour \ $j+k = 0$.

\bigskip

Supposons l'assertion d\'emontr\'ee pour \ $j+k \leq n$ \ et montrons-l\`a pour \ $j+k = n+1$. Supposons, par exemple \ $j \geq 1$. Alors on a, gr\^ace \`a l'hypoth\`ese de r\'ecurrence
$$  G_{\alpha,\beta,j-1,k}^{p,q}(s) - s^p.\bar s^q.\vert s\vert^{2(\alpha+\beta+1)}.P(Log\vert s\vert) \in \mathscr{C}^{\infty}.$$
Comme on a
$$ \frac{\partial G_{\alpha,\beta,j-1.k}^{p,q}}{\partial \alpha} = G_{\alpha,\beta,j,k}^{p,q} $$
et que la d\'ependance en \ $(\alpha,\beta)\in (]-1,+\infty[)^2$ \ des calculs pr\'ec\'edents est \ $\mathscr{C}^{\infty}$, on obtient en d\'erivant en \ $\alpha$ \ la relation \ $(^{**})$ \ pour le couple \ $(j-1,k)$
$$ G_{\alpha,\beta,j,k}^{p,q}(s) -  s^p.\bar s^q.\vert s\vert^{2(\alpha+\beta+1)}.\big[P_{j-1,k}(Log\vert s\vert).Log\vert s\vert^2 + \frac{\partial P_{j-1,k}}{\partial \alpha}(Log\vert s\vert)\big] \in \mathscr{C}^{\infty}$$
ce qui prouve notre assertion pour le couple \ $(j,k)$ \ en vertu de l'unicit\'e du d\'eveloppement asymptotique.

\parag{Remarque} Si \ $\alpha $ \ ou \ $\beta$ \ est entier, mais toujours en supposant que \ $\alpha + \beta +1 \not\in \mathbb{N}$, on constate que pour \ $j + k = 0$ \ on n'a pas de terme singulier pour la fonction \ $G^{p,q}_{\alpha,\beta,0,0}$. La r\'ecurrence sur \ $j+k$ \ montre que l'on peut prendre \ $L = j+k-1$ \ dans ce cas. $\hfill \square$

\bigskip

Dans le cas o\`u \ $\alpha + \beta +1 = m \in \mathbb{N}$ \ la fonction \ $\Phi_{0,0}$ \ sera la somme d'une fonction analytique et d'un terme logarithmique du type \ $ \delta.\vert s\vert^{2m +p+q}.Log \vert s\vert^2 $, qui sera obtenu pour \ $\nu = 2m + p+ q $. On peut alors mettre le terme initial \ $c.\vert s\vert^{2m}.s^p.\bar s^q$ \ dans les termes non singuliers, sortir le terme logarithmique et conclure de fa{\c c}on analogue au cas pr\'ec\'edent pour la fonction \ $G_{\alpha,\beta,0,0}^{p,q}$. Le cas o\`u \ $(j,k)\in \mathbb{N}^2$ \ est arbitraire s'en d\'eduit alors par d\'erivation \ $j-$fois en \ $\alpha$ \ et \ $k-$fois en \ $\beta$.

\bigskip

Il nous reste \`a prouver que pour \ $(\alpha,\beta) \in \mathbb{N}^2$ \ on a en fait un polyn\^ome en \ $Log\vert s\vert$ \ de degr\'e \ $\leq j+k-1$. Il suffit en fait de montrer ce r\'esultat pour \ $j+k = 1$ \ car le cas g\'en\'eral s'en d\'eduit imm\'ediatement par d\'erivation  en \ $\alpha$ \ et \ $\beta$.\\
Il s'agit donc de montrer que pour  \ $(\alpha,\beta) \in \mathbb{N}^2$ \  et, par exemple\footnote{rappelons que l'on a sym\'etrie entre \ $j$ \ et \ $k$.} \ $j=0,k = 1$ \ la fonction \ $G_{\alpha,\beta,0,1}^{p,q}$ \ est \ $\mathscr{C}^{\infty}$. Mais ceci est \'evident puisque la fonction \ $u \to \vert u\vert^{2\alpha}$ \ est \ $\mathscr{C}^{\infty}$ \ dans ce cas. $\hfill \blacksquare$

\bigskip

\parag{Cinqui\`eme et derni\`ere \'etape} D'apr\`es la seconde \'etape, on peut se contenter d'integrer sur le disque \ $\{ \vert u\vert \leq 1\}$ \ pour prouver l'assertion. Fixons un entier \ $N$. D'apr\`es la troisi\`eme \'etape, on peut alors remplacer dans l'int\'egrale les fonctions \ $\theta$ \ et \ $\eta$ \ par des polyn\^omes de degr\'es \ $N +1$ \ pour \'etablir l'existence et la forme du d\'eveloppement asymptotique d'ordre \ $N$ \ puisque l'erreur commise sera de classe \ $\mathscr{C}^N$ \ au voisinage de l'origine.\\
En d\'eveloppant le polyn\^ome en \ $s-u$ \ on est alors ramen\'e \`a montrer l'existence d'un d\'eveloppement asymptotique du type d\'esir\'e pour les fonctions du type 
$$ s^{p'}.\bar s^{q'}.G_{\alpha,\beta,j,k}^{p,q}(s) .$$
Ceci est donn\'e par la proposition de l'\'etape 4.\\
On remarquera que l'on a prouv\'e en m\^eme temps que les fonctions \ $F_{\alpha,\beta,j,k}$ \ sont \ $\mathscr{C}^{\infty}$ \ en dehors de l'origine, puisque c'est le cas pour les fonctions \ $G_{\alpha,\beta,j,k}^{p,q}$. Ceci ach\`eve la preuve du th\'eor\`eme \ref{T-S}. $\hfill \blacksquare$

\subsection{Preuves de la proposition \ref{DA uni}  et du corollaire  \ref{Bern}.}

Commen{\c c}ons par d\'emontrer la proposition \ref{DA uni}.

\parag{Preuve de la proposition \ref{DA uni} } Soit donc \ $K$ \ un compact de \ $X$ \ contenant le support de \ $\varphi$ \ et notons
$$ \mathcal{M}_K : = \{ DA\big[ \int_{f=s} \ \psi ] ,\ \psi \in \mathscr{C}^{\infty}_K(X,\mathbb{C})^{n,n} \} $$
o\`u \ $DA(F)$ \ d\'esigne le d\'eveloppement asymptotique \`a l'origine de la fonction \ $F$. Alors \ $ \mathcal{M}_K$ \ est un \ $\mathbb{C}[[s,\bar s]]-$module de type fini d'apr\`es [B. 82], qui est contenu dans
$$ \vert \Xi\vert^2_{R,n} : = \sum_{(r,j) \in R \times [0,n]} \quad \mathbb{C}[[s,\bar s]].\vert s\vert^{2r}.(Log\vert s\vert )^j ,$$
o\`u \ $R \subset \mathbb{Q} \cap [0,1[$ \ est un sous-ensemble fini.\\
Notons \ $\tilde{\mathcal{M}}_K$ \ le satur\'e de \ $\mathcal{M}_K$ \ par les op\'erateurs \ $s\frac{\partial}{\partial s}$ \ et \ $\bar s\frac{\partial}{\partial \bar s}$. Comme \ $\vert \Xi\vert^2_{R,n}$ \ est stable par ces deux op\'erateurs, la noeth\'erianit\'e de \ $\mathbb{C}[[s,\bar s]]$ \ implique que \ $\tilde{\mathcal{M}}_K$ \ est \'egalement de type fini sur  \ $\mathbb{C}[[s,\bar s]]$.\\
Mais d'apr\`es [B.-S. 74]\footnote{en combinant l'inclusion \ $f^{n+1}.\Omega^{n+1}_X \subset df\wedge \Omega^n_X$ \  avec la formule de d\'erivation  $$\frac{\partial}{\partial s} \int_{f=s} \varphi = \int_{f=s} \frac{d'\varphi}{df} .$$}, on a \ $s^{n+1}.\frac{\partial}{\partial s}  \mathcal{M}_K \subset  \mathcal{M}_K$ \ ainsi que \ $\bar s^{n+1}.\frac{\partial}{\partial \bar s}  \mathcal{M}_K \subset  \mathcal{M}_K$.  Ceci montre que l'on a l'inclusion
$$ \tilde{\mathcal{M}}_K \subset  \mathcal{M}_K[s^{-1}, \bar s^{-1}] .$$
La finitude de \ $ \tilde{\mathcal{M}}_K$ \ sur \ $\mathbb{C}[[s,\bar s]]$ \ donne alors l'existence d'un entier \ $\kappa$, ne d\'ependant que de \ $f$ \ et du compact \ $K$,  tel que l'on ait 
$$ \vert s\vert^{2\kappa}. \tilde{\mathcal{M}}_K \subset \mathcal{M}_K .$$
Il nous suffit donc de prouver la proposition pour \ $\tilde{\mathcal{M}}_K$ \ en prenant \ $\kappa = 0$ \ dans ce cas.\\

Comme \ $\tilde{\mathcal{M}}_K$ \ est stable par \ $s\frac{\partial}{\partial s}$ \ et \ $\bar s\frac{\partial}{\partial \bar s}$, l'utilisation it\'er\'ee des op\'erateurs
$$ s\frac{\partial}{\partial s} - (r_1+m_1) \quad {\rm et} \quad \bar s\frac{\partial}{\partial \bar s} - (r_1+m'_1)$$
permet de supprimer dans le d\'eveloppement asymptotique de \ $\varphi$ \ tous les termes qui ne sont pas des \ $\mathcal{O}(\vert s\vert^{2N})$  (il n'y en a qu'un nombre fini) et qui ne sont pas de la forme \ $\vert s\vert^{2r}.s^m.\bar s^{m'}.P(Log\vert s\vert)$ \ o\`u \ $P$ \ est un polyn\^ome de degr\'e \'egal \`a \ $j$ \ (rappelons que l'on a suppos\'e que \ $j$ \ est maximal pour \ $(r,m,m')$ \ donn\'e). Ceci permet alors de conclure. $\hfill \blacksquare$

\parag{Remarque} Si la fonction \ $f$ \ est contenue dans son id\'eal jacobien au voisinage du compact \ $K$, on a \ $\tilde{\mathcal{M}}_K = \mathcal{M}_K$ \ et l'on peut prendre \ $\kappa = 0$ \ dans la proposition. On a donc le r\'esultat optimal dans ce cas.

\bigskip

En utilisant de mani\`ere anticip\'ee les r\'esultats de la section suivante, nous allons d\'emontrer le  le corollaire \ref{Bern}. 

\parag{Preuve du corollaire \ref{Bern}} La premi\`ere partie du corollaire est imm\'ediate puisque l'on sait que modulo \ $\mathbb{Z}$ \ les racines du polyn\^ome de Bernstein de \ $f$ \ \`a l'origine correspondent aux exposants qui apparaissent effectivement dans les\\
 d\'eveloppements asymptotiques des int\'egrales fibres au voisinage de l'origine ; voir par exemple [Bj. 93]  chapitre 6.5.\\
La r\'eciproque est cons\'equence  de la proposition \ref{DA uni}  qui permet  d'appliquer le th\'eor\`eme pr\'ecis \ref{Precis} et ainsi d'\'eviter les ph\'enom\`enes de compensation entre diff\'erents termes des d\'eveloppements asymptotiques de \ $f$ \ et \ $g$. $\hfill \blacksquare$

\section{Calcul des constantes : \\
d\'emonstration du th\'eor\`eme pr\'ecis.}

\subsection{Pr\'eliminaires.}

  \parag{Notation} Pour \ $a, b, \cdots$ \ des nombres complexes, nous noterons \ $\alpha, \beta, \cdots$ \ leurs parties r\'eelles respectives.\\
  
  Notons pour \ $p,q$ \ entiers positifs,
      \begin{align*}
      & U_{p,q} : = \{ (a,b)\in \mathbb{C}^2 \ / \  \alpha + p/2 > -1 \quad {\rm et} \quad \beta + q/2 > -1 \} \\
      & U^0_{p,q} : = U_{p,q} \cap \{(a,b)\in \mathbb{C}^2 \ / \ \alpha + \beta +1+ \frac{p+q}{2} < 0 \}\\
      & V_{p,q} : = U_{p,q} \setminus \{(a,b)\in \mathbb{C}^2 \ / \ a+b+1 \in \mathbb{Z} \}.
      \end{align*}
      On remarquera que pour \ $(a,b) \in U^0_{p,q}$, on a \ $a+b+1 \not\in \mathbb{N}$ \ puisque les in\'egalit\'es impos\'ees impliquent \ $-1-\frac{p+q}{2} < \alpha+\beta+1 < -\frac{p+q}{2}$. $\hfill \square$
      
      \bigskip

 Nous poserons, pour \ $ x \in \mathbb{R}, \vert x\vert  < 1, m \in \mathbb{Z}$ \ et \ $a \in \mathbb{C}$ \ tel que \ $\alpha > -1$
       \begin{equation*}
        \frac{1}{2\pi}\int_0^{2\pi} \vert 1 - x.e^{-i\theta}\vert^{2a}.e^{im\theta}.d\theta = \sum_{r=0}^{+\infty} \gamma_{a,m}^r.x^{r} . \tag{@}
        \end{equation*}
       On remarquera que \ $  \gamma_{a,m}^r = 0 $ \ si \ $r \not= m $ \ modulo 2, puisque l'int\'egrale consid\'er\'ee est invariante par le changement \ $x \to - x$ \ et \ $\theta \to \theta + \pi$. \\
       La conjuguaison complexe montre que l'on a
       $$ \overline{\gamma_{\bar a,-m}^r} = \gamma_{a,m}^r $$
       il nous suffira donc de consid\'erer le cas o\`u \ $m$ \ est dans \ $\mathbb{N}$.\\
       On remarquera \'egalement que les coefficients \ $ \gamma_{a,m}^r $ \ d\'ependent holomorphiquement de \ $a$ \ sur l'ouvert \ $\alpha > -1$, d'apr\`es la formule de Cauchy.
       
         \begin{lemma}       
        La fonction holomorphe d\'efinie sur l'ouvert  \ $U_q^0 : = U^0_{0,q}$ \ de \ $\mathbb{C}^2$ \  en posant
       $$ G_{q}(a,b) : = \frac{i}{4\pi} \int_{\mathbb{C}} \vert 1-t\vert^{2a}.t^q.\vert t\vert^{2b}.dt\wedge d\bar t $$
       est \'egale \`a 
       \begin{equation*} \frac{1}{2}.\frac{\Gamma(a+1)\Gamma(b+q+1)\Gamma(-a-b -1)}{\Gamma(-a)\Gamma(-b)\Gamma(a+b+q+2)}. \tag{*}
       \end{equation*}
       \end{lemma}
       
       \parag{Preuve} Commen{\c c}ons par remarquer que la convergence absolue de l'int\'egrale d\'efinissant \ $G_{q}$ \ est assur\'ee  en \ $t = 0$ \ par la condition \ $2\beta + q +1  > -1$, en \ $t = 1$ \ par la condition \ $\alpha > -1$ \ et en \ $\vert t\vert = + \infty $ \ par la condition \ $ 2(\alpha + \beta) + q + 1 < -1$. Donc la fonction \ $G_{q}$ \ est bien holomorphe sur l'ouvert \ $U^0_q$. \\
       Pour prouver la formule \ $(^*)$ \ il suffit de le faire quand \ $a = \alpha$ \ et \ $b = \beta$. Ce que nous supposerons donc dans la suite.\\
       En coordonn\'ees polaires \ $ t = \rho.e^{i\theta}$ \ on aura, pour \ $(\alpha,\beta) \in U^0_q$ :
       \begin{align*}
       & G_{q}(\alpha,\beta) = \int_0^{+\infty} (1+\rho^2)^{\alpha}.\rho^{2(\beta+q/2+1)}.\frac{d\rho}{\rho}.\Big(\frac{1}{2\pi}\int_0^{2\pi} \big(1 -\frac{2\rho.cos\,\theta}{1+\rho^2}\big)^{\alpha}.e^{iq\theta}.d\theta \Big) \tag{**}
       \end{align*}
       Mais on a
       \begin{align*}
       & \frac{1}{2\pi}\int_0^{2\pi} \big(1 -\frac{2\rho.cos\,\theta}{1+\rho^2}\big)^{\alpha}.e^{iq\theta}.d\theta = \\
       & \qquad \sum_{k=0}^{\infty} \frac{\Gamma(\alpha+1)}{\Gamma(\alpha-k+1).k!}.\frac{(-1)^k.2^k.\rho^k}{(1+\rho^2)^k}.\frac{1}{2\pi}\int_0^{2\pi} cos^k\theta.e^{iq\theta}.d\theta .
       \end{align*}
       Et, comme on a
       \begin{align*}
       & \frac{1}{2\pi}\int_0^{2\pi} cos^k\theta.e^{iq\theta}.d\theta \ 
= \frac{1}{2^k}.C_k^j \quad \quad {\rm si } \quad k = q + 2j ,  \ {\rm avec} \quad j \in [0,k], \ {\rm et} \\
& \qquad  \qquad \qquad \qquad \quad \ = 0 \quad {\rm sinon, \ on \ obtient}
\end{align*}

\begin{align*}
 & \frac{1}{2\pi}\int_0^{2\pi} \big(1 -\frac{2\rho.cos\,\theta}{1+\rho^2}\big)^{\alpha}.e^{iq\theta}.d\theta  = \\
 & \qquad \qquad  (-1)^q.\sum_{j=0}^{+\infty} C_{q+2j}^j\times\frac{\Gamma(\alpha+1)}{\Gamma(\alpha-q -2j +1).(q+2j)!}\times\frac{\rho^{q+2j}}{(1+\rho^2)^{q+2j}}.
 \end{align*}
 En reportant dans \ $(^{**})$ \ cela donne
 \begin{align*}
 & G_{q}(\alpha,\beta) = \\
 &  (-1)^q.\sum_{j=0}^{+\infty} C_{q+2j}^j\times \frac{\Gamma(\alpha+1)}{\Gamma(\alpha-q -2j +1).(q+2j)!}\times\int_0^{+\infty} (1+\rho^2)^{\alpha -q-2j}.\rho^{2(\beta+q+j+1)}.\frac{d\rho}{\rho}.
  \end{align*}
  Utilisons alors la formule classique
  $$ \int_0^{+\infty} (1+x^2)^{-u}.x^v.dx = \frac{1}{2}.\frac{\Gamma(\frac{v+1}{2}).\Gamma(u -\frac{v+1}{2})}{\Gamma(u)} $$
  avec \ $ u = q + 2j - \alpha$ \ et \ $1+v = 2(\beta+q+j+1) $. On obtient
  \begin{align*}
  & G_{q}(\alpha,\beta) = \frac{(-1)^q}{2}.\Gamma(\alpha+1).\sum_{j=0}^{+\infty} \ \frac{\Gamma(\beta+q+j+1).\Gamma(-\alpha-\beta+j-1).(q+2j)!}{\Gamma(\alpha-q-2j+1).[(q+2j)!]\Gamma(q+2j-\alpha).j! (q+j)!}.
  \end{align*}
  Apr\`es simplification par \ $(q+2j)!$ \ et utilisation de la formule des compl\'ements qui donne
  \begin{align*}
  & \Gamma(\alpha-q-2j+1).\Gamma(q+2j-\alpha)  = \frac{\pi}{\sin \pi(q+2j-\alpha)} \\
  & \qquad  = (-1)^{q+1} \frac{\pi}{\sin\pi\alpha} = (-1)^{q}.\Gamma(1+\alpha).\Gamma(-\alpha)
  \end{align*}
  on arrive \`a 
  \begin{align*}
   & G_{q}(\alpha,\beta) =  \frac{1}{2.\Gamma(-\alpha)}\times \sum_{j=0}^{+\infty} \ \frac{\Gamma(\beta+q+j+1).\Gamma(-\alpha-\beta+j-1)}{\Gamma(q+j+1).j!}
   \end{align*}
   Utilisons maintenant la formule (voir Erdelyi, Magnus ... p. 10)
   $$ \sum_{j=0}^{+\infty} \ \frac{\Gamma(j+x).\Gamma(j+y)}{\Gamma(j+z).j!} = \frac{\Gamma(x).\Gamma(y).\Gamma(z-x-y)}{\Gamma(z-x).\Gamma(z-y)} .$$
   On obtient avec \ $x = \beta+q+1, y = -\alpha-\beta-1$ \ et \ $z = q+1$
   $$  G_{q}(\alpha,\beta) \  = \frac{1}{2} . \frac{\Gamma(\beta+q+1).\Gamma(-\alpha-\beta-1).\Gamma(\alpha+1)}{\Gamma(-\alpha).\Gamma(-\beta).\Gamma(\alpha+\beta+q+2)}$$
   c'est \`a dire la formule annonc\'ee.$\hfill \blacksquare$
   
   \parag{Remarques}
   \begin{enumerate}[1)]
   \item Le changement \ $t \to 1/t$ \ transforme \ $G_{q}(\alpha,\beta)$ \ en \ $\overline{G_{q}(\alpha,-(\alpha+\beta+q+2))}$. \\
   On constate que ceci est bien compatible avec la formule obtenue puisque si \\
    $\alpha' = \alpha, \ \beta' = -(\alpha+\beta+q+2)$ \ on a
   \begin{align*}
   & \beta' + q + 1 = -\alpha-\beta-1, \quad -\alpha'-\beta'-1 = \beta+q+1,\\
   &   \alpha'+\beta'+q+2 = -\beta  \quad \quad   etc...
    \end{align*}
   \item On notera \'egalement que l'ouvert \ $U^0_{0,q}$ \ est stable par cette involution, puisque \\
     $$\beta' + q/2 +1= -\alpha-\beta -q/2- 1 > 0  \quad {\rm ainsi \ que}  $$
     $$\alpha'+\beta'+q/2 +1 = -\beta -q/2 -1 < 0 $$
      pour \ $(\alpha,\beta) \in U^0_{0,q}$. 
      \item On remarquera que pour \ $(\alpha,\beta) \in U^0_{0,q} \cap  \mathbb{R}^2$ \ le nombre \ $G_{q}(\alpha,\beta)$ \ est r\'eel. On en d\'eduit que l'on a, pour tout \ $(\alpha,\beta) \in U^0_{0,q} \cap  \mathbb{R}^2$ 
      \begin{equation*}
       G_{\bar q}(\alpha,\beta) : = \frac{i}{4\pi} \int\int_{\mathbb{C}} \vert 1-t\vert^{2\alpha}.\bar{t}^q.\vert t\vert^{2\beta}.dt\wedge d\bar t  = G_{q}(\alpha,\beta).
      \end{equation*}
      On a donc, pour \ $(a,b) \in U^0_{0,q}$ \ l'\'egalit\'e \  $   G_{\bar q}(a,b) = \overline{G_{q}(\bar a,\bar b)}$. $\hfill \square $
      \end{enumerate}
      
      Consid\'erons maintenant la m\^eme int\'egrale que dans le lemme pr\'ec\'edent mais en demandant seulement \`a \ $(a,b)$ \ de v\'erifier les conditions de convergence en \ $t = 0$ \ et \ $t =1$ \ c'est \`a dire d'\^etre dans l'ouvert \ $U_{0,q}$ \ d\'efini par les in\'egalit\'es \ $\alpha > -1$ \ et \ $\beta + q/2 > -1$. L'int\'egrale diverge alors \`a l'infini, mais le lemme suivant montre que sa "partie finie", c'est \`a dire le terme constant dans le d\'eveloppement asymptotique \`a l'infini de cette int\'egrale donne le prolongement m\'eromorphe de la fonction \ $G_q$ \ \`a l'ouvert \ $U_{0,q}$.
      
       \begin{lemma}\label{Calcul de la constante}
      Soit  \ $q$ \ un entier positif ou nul. Pour \ $(a,b) \in V_{0,q}$ \  on a, pour \ $N$ \ entier assez grand, 
      \begin{align*}
      & \tag{***}  \lim_{s \to 0} \Big[ \frac{i}{4\pi}\int_{\vert t \vert \leq 3} \vert 1-t\vert^{2a}.t^{q}.\vert t\vert^{2b}.dt\wedge d\bar t  \ + \\
      & \qquad \int_3^{1/\vert s\vert} \rho^{2(a+b+1)+q}.\frac{d\rho}{\rho}.\big[\frac{1}{2\pi}\int_0^{2\pi} \vert1-\frac{e^{-i\theta}}{\rho}\vert^{2a}.e^{iq\theta}.d\theta - \sum_0^N \gamma_{a,q}^r.\rho^{-r}\big]\Big] \  =  \\
      &   \frac{1}{2}. \frac{\Gamma(a+1).\Gamma(b+q+1).\Gamma(-a-b-1)}{\Gamma(-a).\Gamma(-b).\Gamma(a+b+q+2)}  +  \sum_{r=0}^N \frac{\gamma_{a,q}^r}{2(a+b+1)+q-r}.3^{2(a+b+1)+q-r}  
      \end{align*}
      \end{lemma}

      \parag{Preuve} Remarquons d\'ej\`a que pour \ $N$\, fix\'e, la diff\'erence
      $$\big[\frac{1}{2\pi}\int_0^{2\pi} \vert1-\frac{e^{-i\theta}}{\rho}\vert^{2a}.e^{iq\theta}.d\theta - \sum_0^N \gamma_{a,q}^r.\rho^{-r}\big] $$
      est un \ $O(\rho^{-N-1})$ \ et donc que pour \ $2(\alpha+\beta+1) + q-N-1 < 0$ \ l'int\'egrale de \ $3$  \ \`a \ $+\infty$ \ converge et la limite cherch\'ee est simplement la valeur de cette int\'egrale.\\
      Il s'agit donc simplement de montrer que le prolongement analytique de la fonction holomorphe \ $G_{q}$ \ d\'efinie sur l'ouvert \ $U^0_{0,q}$ \ est bien donn\'e par cette int\'egrale. On concluera alors  gr\^ace au lemme pr\'ec\'edent.\\
       Mais pour \ $(a,b) \in U^0_{0,q}$, en coupant l'int\'egrale d\'efinissant \ $G_q$ \ pour \ $\vert t\vert \leq 3$ \ et pour \ $\vert t \vert \geq 3$ \ et en utilisant le d\'eveloppement \ $(@)$ \ \`a l'ordre \ $N$ \ dans cette derni\`ere, on obtient bien l'\'egalit\'e de la limite du membre de gauche de \ $(^{***})$ \ avec \ $G_q$ \ sur l'ouvert \ $U^0_{0,q}$, puisque \ $\vert s\vert^{-[2(a+b+1)+q-r]}$ \ tend vers \ $0$ \ quand \ $s \to 0$, \'etant donn\'e que l'on a \ $\alpha+ \beta+q/2 + 1 < 0$ \ et \ $r\geq 0$. $\hfill \blacksquare$
       
       \bigskip
       
       On remarquera que ceci montre que le prolongement analytique de la fonction  \ $G_{q}$ \ a des p\^oles  simples (au plus) sur les droites 
        $$ a+b+1 = \frac{r - q}{2} \quad {\rm pour} \quad  r \in \mathbb{N}, \quad r = q \quad modulo \ 2.$$
        Le fait que ces points soient r\'eellement des p\^oles simples quand \ $a+b+1 \in \mathbb{N}$ \ est montr\'e \`a posteriori par la formule que l'on a \'etablie. Ceci permet de calculer les coefficients \ $\gamma_{a,q}^r $. On obtient, pour \ $a \not\in \mathbb{N}$, en utilisant la formule des compl\'ements (on notera  que \ $a \not\in -\mathbb{N}$ \ sous nos hypoth\`eses, puisque l'on a \ $\alpha > -1$ \ et aussi que \ $\frac{r+q}{2} - a \not\in -\mathbb{N}$, puisque \ $\frac{r+q}{2} - \alpha = \beta+q+1 > 0$)
        \begin{align*}
        &  \gamma_{a,q}^r =  \frac{\Gamma(a+1).\Gamma(\frac{r+q}{2}-a).(-1)^{\frac{r-q}{2}}}{\Gamma(-a).\Gamma(a+1-\frac{r-q}{2}).\Gamma(\frac{r+q}{2}+1).(\frac{r-q}{2})!} \\
        & \qquad = (-1)^r.\frac{\Gamma(a+1)^2}{\Gamma(a+1-\frac{r-q}{2}).\Gamma(a+1-\frac{r+q}{2}).(\frac{r-q}{2})! (\frac{r+q}{2})!}\quad  . \tag{@@}
        \end{align*}
         Pour \ $a \in \mathbb{N}$ \   le calcul direct donne, pour \ $q \leq r \leq r+q \leq 2a$ 
        $$  \gamma_{a,q}^r = (-1)^r.C_a^{(r+q)/2}.C_a^{(r-q)/2},$$  
        et \ $0$ \ sinon, ce qui co{\"i}ncide bien avec \ $(@@)$.  
        
          \bigskip

         La situation pr\'ec\'edente  o\`u \ $a+b+1 \in \mathbb{N}$ \ avec \ $a$ \ non entier  est examin\'ee dans le lemme suivant.
         
         \bigskip
               
        \begin{lemma}\label{Calcul de la constante bis}
         Soit  \ $q$ \ un entier positif ou nul. Pour \ $(a,b) \in U_{0,q}, a+b+1 \in \mathbb{N}$, \ $a$ \ {\bf non entier},  posons \ $r_0 = 2(a+b+1) + q $. Alors on a, pour \ $N \geq r_0$ \ entier    
            \begin{align*}
      & \lim_{s \to 0} \Big[ \frac{i}{4\pi}\int_{\vert t \vert \leq 3} \vert 1-t\vert^{2a}.t^{q}.\vert t\vert^{2b}.dt\wedge d\bar t  \ + \\
      & \qquad \int_3^{1/\vert s\vert} \rho^{2(a+b+1)+q+j}.\frac{d\rho}{\rho}.\big[\frac{1}{2\pi}\int_0^{2\pi} \vert1-\frac{e^{-i.\theta}}{\rho}\vert^{2a}.e^{iq\theta}.d\theta - \sum_0^N \gamma_{a,q}^r.\rho^{-r}\big]\Big] \  + \\
      &\qquad \qquad - \sum_{r=0, r \not= r_0}^N \frac{\gamma_{a,q}^r}{2(a+b+1)+q-r}.3^{2(a+b+1)+q-r} + \gamma_{a,q}^{r_0}.Log\, 3 
      \end{align*}
      est \'egale \`a          
       $$ \frac{(-1)^{r_0}}{2}\frac{\big[\Gamma'(1) + \sum_{j=1}^{(r_0 -q)/2} \frac{1}{j}\big].\Gamma(a+1)^2}{\Gamma((a+1-\frac{r_0+q}{2}).\Gamma((a+1-\frac{r_0-q}{2}).(\frac{r_0+q}{2}!).(\frac{r_0-q}{2}!)} \quad .$$
         \end{lemma}

      \parag{Preuve} La preuve est la m\^eme que pour le lemme pr\'ec\'edent sauf que le terme pour \ $r = r_0$ \ donne un terme logarithmique dans ce cas. On doit donc \'evaluer la valeur en \ $2(a+b+1)+q = r_0$ \ de la somme
      $$\frac{1}{2}\, \frac{\Gamma(a+1).\Gamma(b+q+1).\Gamma(-a-b-1)}{\Gamma(-a).\Gamma(-b).\Gamma(a+b+q+2)}  \ + \   \frac{1}{2}\,\frac{\gamma_{a,q}^{r_0}}{a+b+1+(q-r_0)/2}. $$
                
         Compte tenu de la formule \ $(@@)$,  du fait que l'on a, pour \ $k \in \mathbb{N}$,
         $$ \Gamma(z-k) = \frac{(-1)^k}{k!}.[\frac{1}{z} + \Gamma'(1) + \sum_{j=1}^k \frac{1}{j} + o(z)] $$
         quand \ $ z \to 0$, un calcul simple montre que la valeur de cette limite est bien celle annonc\'ee.
 $\hfill \blacksquare$
      
      \bigskip
      
      On notera que cette limite n'est jamais nulle puisque l'on suppose que \ $a$ \ n'est pas un entier

        \bigskip
        
         \begin{lemma}\label{Calcul des constantes} 
      Soient \ $0 \leq p \leq q$ \ deux entiers. Posons pour \ $(a,b)\in U^0_{p,q}$
      \begin{align*}
       &  F_{p,q}(a,b) : = \frac{i}{\pi} \int\int_{\mathbb{C}} \vert 1-t\vert^{2a}.(1-t)^p.\vert t\vert^{2b}.t^q.dt\wedge d\bar t  \\
  &  F_{p,\bar q}(a,b) : = \frac{i}{\pi} \int\int_{\mathbb{C}} \vert 1-t\vert^{2a}.(1-t)^p.\vert t\vert^{2b}.\bar t^q.dt\wedge d\bar t  
      \end{align*}
           On a  alors
      $$ F_{p,q}(a,b) = \frac{\Gamma(a+p+1).\Gamma(b+q+1).\Gamma(-a-b-1)}{\Gamma(a+b+p+q+2).\Gamma(-a).\Gamma(-b)} $$
      $$  F_{p,\bar q}(a,b) = (-1)^p .\frac{\Gamma(a+p+1).\Gamma(b+q+1).\Gamma(-a-b-p-1)}{\Gamma(a+b+q+2).\Gamma(-a).\Gamma(-b)} \quad . $$
      \end{lemma}
      
      \parag{Remarques} 
      \begin{enumerate}[i)]
      \item  On v\'erifie facilement  les \'egalit\'es "\'evidentes \`a priori" sur \ $U^0_{p,q} \cap \mathbb{R}^2$
         $$ F_{p,0} = F_{p,\bar 0} \quad {\rm et} \quad F_{0,q} = F_{0,\bar q} . $$
         \item Une cons\'equence simple du lemme pr\'ec\'edent est le prolongement m\'eromorphe des fonctions \ $F_{p,q}$ \ et \ $F_{p,\bar q}$ \ \`a l'ouvert \ $U_{p,q}$ \ avec des p\^oles au plus simples sur les droites \ $a + b+ 1 \in \mathbb{N}$.
         \end{enumerate}

      \parag{Preuve} La formule suivante sera la clef de cette preuve. 
      
      \begin{lemma}\label{Formule}
      Soit \ $p \in \mathbb{N}$. Pour \ $x,y$ \ dans \ $\mathbb{C}\setminus \{-\mathbb{N}\}$ \ on a
      $$ A_p(x,y) : = \sum_{j=0}^p (-1)^j.C_p^j.\frac{\Gamma(x+j)}{\Gamma(x+y+j)} = \frac{\Gamma(x).\Gamma(y+p)}{\Gamma(x+y+p).\Gamma(y)}\quad .$$
      \end{lemma}
      
      \parag{Preuve} Montrons cette \'egalit\'e par r\'ecurrence sur \ $p \geq 0$. Le cas \ $p = 0$ \ est trivial. Supposons la formule d\'emontr\'ee pour \ $p$ \ et montrons-l\`a pour \ $p+1$. En utilisant l'\'egalit\'e \ $C_{p+1}^j = C_p^j + C_p^{j-1} $ \ on obtient
      \begin{align*}
      & A_{p+1}(x,y) = A_p(x,y) - A_p(x+1,y) = \frac{\Gamma(x).\Gamma(y+p)}{\Gamma(x+y+p).\Gamma(y)} - \frac{\Gamma(x+1).\Gamma(y+p)}{\Gamma(x+y+p+1).\Gamma(y)}\\
      & \quad \quad = \frac{\Gamma(x).\Gamma(y+p)}{\Gamma(x+y+p+1).\Gamma(y)}[x+y+p - x] \\
       & \quad \quad =  \frac{\Gamma(x).\Gamma(y+p+1)}{\Gamma(x+y+p+1).\Gamma(y)}\quad .  \qquad \qquad \qquad\qquad \qquad \qquad \qquad \qquad \qquad\qquad \blacksquare 
       \end{align*}
       
      \parag{Preuve du lemme \ref{Formule}} La formule du bin\^ome donne
       \begin{align*}
       &  F_{p,q}(a,b) = \sum_{j=0}^p (-1)^j.C_p^j.G_{q+j}(a,b) \\
       & \qquad \qquad = \frac{\Gamma(a+1).\Gamma(-a-b-1)}{\Gamma(-a).\Gamma(-b)}\times \Big(\sum_{j=0}^p (-1)^j.C_p^j.\frac{\Gamma(b+q+j+1)}{\Gamma(a+b+q+j+2)} \Big)
       \end{align*}
       et en utilisant le lemme ci-dessus avec \ $x = b+q+1, y= a+1 $ \ on obtient la formule annonc\'ee.
       
       \smallskip
       
       La formule du bin\^ome donne
       \begin{align*}
       & F_{p,\bar q}(a,b) = \sum_{j=0}^p (-1)^j.C_p^j.G_{q-j}(a,b+j) \\
       & \qquad \qquad = \frac{\Gamma(a+1).\Gamma(b+q+1)}{\Gamma(-a).\Gamma(a+b+q+2)}\times\Big(\sum_{j=0}^p (-1)^j.C_p^j.\frac{\Gamma(-a-b-j-1)}{\Gamma(-b-j)} \Big)
       \end{align*}
       La formule des compl\'ements donne alors la formule annonc\'ee. \ $\hfill \blacksquare$
       
       \parag{Remarque} Pour \ $a+b+1 \not\in \mathbb{N}$, les nombres complexes \ $F_{p,q}(a,b)$ \ et \ $F_{p,\bar q}(a,b)$ \ sont non nuls. $\hfill \square$
       
\bigskip

\subsection{Le premier cas.}

\parag{Notation} Soient \ $0 \leq p \leq q$ \ deux entiers  et soit \ $D$ \ le disque unit\'e du plan complexe.  Pour \ $(a,b) \in U_{p,q}$ \ et \ $s \in D^* :  =  D \setminus \{0\}$ \ posons
       \begin{align*}
       & F_{p,q}(a,b)[s] : = \frac{i}{4\pi}\int_{\vert u \vert \leq 1} \vert s-u\vert^{2a}.(s-u)^p.\vert u \vert^{2b}.u^q.du \wedge d\bar u \\
       & F_{p,\bar q}(a,b)[s] : =  \frac{i}{4\pi}\int_{\vert u \vert \leq 1} \vert s-u\vert^{2a}.(s-u)^p.\vert u \vert^{2b}.\bar u^q.du \wedge d\bar u
       \end{align*}

        \begin{prop}[Le premier cas :  \ $a+b+1$ \ n'est pas dans \ $\mathbb{N}$.]\label{Cas 1}
      On suppose que  \ $(a,b) \in V_{p,q}$. Alors  il existe des fonctions \ $\Phi_{p,q}$ \ et \ $\Phi_{p,\bar q}$ qui sont \ $\mathscr{C}^{\infty}$ \ sur \ $V_{p,q} \times D$, holomorphes sur \ $V_{p,q}$ \ pour \ $s \in D$ \ fix\'e, telles que l'on ait sur \ $V_{p,q}\times D^*$ 
       \begin{align*}
       & F_{p,q}(a,b)[s] = F_{p,q}(a,b).s^{p+q}.\vert s\vert^{2(a+b+1)} + s^{p+q}.\Phi_{p,q}(a,b,s) \\
       &  F_{p,\bar q}(a,b)[s] = F_{p,\bar q}(a,b).s^p.\bar s^q.\vert s\vert^{2(a+b+1)} + \bar s^{q-p}.\Phi_{p,\bar q}(a,b,s) 
       \end{align*}
       \end{prop}
       
       \parag{Preuve} Remarquons d\'ej\`a que le fait que ces fonctions soient \ $\mathscr{C}^{\infty}$ \ sur l'ouvert  \ $U^0_{p,q}\times D^*$, holomorphes sur \ $U^0_{p,q}$ \ pour \ $s \in D^*$ \ fix\'e, est cons\'equence imm\'ediate des d\'efinitions. Nous allons commencer par traiter le cas de \ $F_{p,q}$. Pour \ $s \in D^*$ \ effectuons le changement de variable \ $u = s.t$. On obtient
       \begin{align*}
       &    F_{p,q}(a,b)[s]  = s^{p+q}.\vert s\vert^{2(a+b+1)}.\frac{i}{4\pi}\int_{\vert t \vert \leq 1/\vert s\vert} \vert 1-t\vert^{2a}.(1-t)^p.\vert t\vert^{2b}.t^q.dt\wedge d\bar t \\
       & \qquad = s^{p+q}.\vert s\vert^{2(a+b+1)}. \sum_{j=0}^p (-1)^j.C_p^j.\frac{i}{4\pi}\int_{\vert t \vert \leq 1/\vert s\vert} \vert 1-t\vert^{2a}.t^{q+j}.\vert t\vert^{2b}.dt\wedge d\bar t .
       \end{align*} 
       Supposons \ $\vert s \vert  < 1/3$ \ et posons 
       $$ C^0_{p,q}(a,b): = \frac{i}{4\pi}\int_{\vert t \vert \leq 3} \vert 1-t\vert^{2a}.(1-t)^p.\vert t\vert^{2b}.t^q.dt\wedge d\bar t .$$
       L'int\'egrale  pour \ $3 \leq \vert t \vert \leq 1/\vert s\vert $ \ correspondant au terme \ $j \in [0,p]$ \ donne en coordonn\'ees polaires
       \begin{equation*}
       I_j =  \int_3^{1/\vert s\vert} \rho^{2(a+b)+q+j+1}.d\rho \frac{1}{2\pi}\int_0^{2\pi} \vert 1-\frac{e^{-i\theta}}{\rho}\vert^{2a}.e^{i(q+j)\theta}.d\theta.  \tag{A}
       \end{equation*}
       Utilisons le d\'eveloppement \ $(@)$, en se souvenant que l'on a \ $r = q+j$ \ modulo 2 :
       \begin{align*}
       & I_j = \sum_{r=0}^{+\infty} \gamma_{a,q+j}^r.\int_0^{1/\vert s\vert} \rho^{2(a+b+(q+j-r)/2))+1}.d\rho \\
       & \quad = \sum_{r=0}^{+\infty} \frac{\gamma_{a,q+j}^r}{2(a+b+(q+j-r)/2)+1)}.\Big[\vert s\vert ^{-2(a+b+(q+j-r)/2)+1)} - 3^{-2(a+b+(q+j-r)/2)+1)}\Big]  \\
       & \quad = C^1_{j,q}(a,b) + \vert s\vert^{-2(a+b+1)}.\Psi_j(a,b,\vert s\vert^2) .  \tag{B} 
       \end{align*}
       la fonction \ $\Psi_j$ \ \'etant une s\'erie de Laurent en \ $\vert s\vert^2$ \ dont les coefficients sont holomorphes en \ $(a,b) $ \ pourvu que \ $a+b \not\in \mathbb{Z}$. On notera que la puissance maximale n\'egative en \ $\vert s\vert^2$ \ dans \ $\Psi_j$ \ est \ $(j+q)/2 \leq (p+q)/2$.\\
       On obtient alors
       $$ F_{p,q}(a,b)[s] =  C_{p,q}(a,b).s^{p+q}.\vert s\vert^{2(a+b+1)} + \sum_{j=0}^p (-1)^j.C_p^j.s^{p+q}.\Psi_j(a,b,\vert s\vert^2)$$
       o\`u nous avons pos\'e
       $$ C_{p,q}(a,b) = C^0_{p,q}(a,b) +  \sum_{j=0}^p (-1)^j.C_p^j.C^1_{j,q}(a,b) .$$
       Ceci \'etablit l'assertion annonc\'ee pour la fonction \ $F_{p,q}$ \ sauf qu'il nous reste \`a montrer les deux points suivants :
       \begin{enumerate}
       \item La fonction \ $\Phi_{p,q}(a,b,s) : = \sum_{j=0}^p (-1)^j.C_p^j.s^{p+q}.\Psi_j(a,b,\vert s\vert^2)$ \ est bien \ $\mathscr{C}^{\infty}$ \ en \ $s=0$, c'est \`a dire ne pr\'esente pas de termes non nul de la forme \ $s^{p+q}/\vert s\vert^{2k}$ \ avec \ $k\geq 1$.
       \item Montrer que la constante \ $C_{p,q}(a,b)$ \ est bien \'egale \`a \ $F_{p,q}(a,b)$.
       \end{enumerate}
       
       Pour \'etablir le premier point consid\'erons un nombre complexe \ $\lambda \in \mathbb{C}^*$ \ et calculons la diff\'erence \ $F_{p,q}(a,b)[\lambda.s] - \lambda^{p+q}.\vert \lambda\vert^{2(a+b+1)}.F_{p,q}(a,b)[s] $. Le changement de variable \ $v = \lambda.u$ \ montre que cette diff\'erence est donn\'ee, pour \ $\vert \lambda \vert < 1$ \  par l'int\'egrale
       $$ \frac{i}{4\pi}\int_{1 \leq \vert v\vert \leq 1/\vert\lambda\vert} \vert s-v\vert^{2a}.(s-v)^p.\vert v\vert^{2b}.v^q.dv\wedge d\bar v .$$
       Cette diff\'erence est donc \ $\mathscr{C}^{\infty}$ \ sur \ $\mathbb{C}^2\times D$, holomorphe sur \ $\mathbb{C}^2$ \  \`a \ $s \in D$ \ fix\'e. Mais la pr\'esence d'un terme non nul de la forme \ $c.s^{p+q}/\vert s\vert^{2k}$ \ avec \ $k \geq 1$ \ produirait dans le d\'eveloppement \`a l'origine de la diff\'erence pr\'ec\'edente le terme 
        $$c.\Big[ \lambda^{p+q}.(\vert \lambda\vert^{-2k} - \vert \lambda\vert^{2(a+b+1)})\Big].s^{p+q}/\vert s\vert^{2k} $$
   n\'ecessairement non nul  puisque \ $\alpha+\beta+1 > -1$, ce  qui contredirait l'aspect \ $\mathscr{C}^{\infty}$ \ en \ $s = 0$ \ de cette diff\'erence.
        
        \bigskip
        
        L'identification de la constante s'obtient facilement en utilisant les lemmes \ref{Calcul de la constante} et \ref{Calcul des constantes}.\\
        La preuve pour la fonction \ $F_{p,\bar q}$ \ est tout \`a fait analogue.\\
        On notera que le terme \ $s^p.\bar s^q/\vert s \vert^{2k}$ \ est \ $\mathscr{C}^{\infty}$ \ pour \ $k \in [0,p]$, puisque l'on a suppos\'e \ $0 \leq p \leq q$, ce qui explique le facteur \ $\bar s^{q-p}$ \ dans ce cas. $\hfill \blacksquare$
        
        \bigskip
        
       \parag{Remarque} Pour \ $a$ \ ou \ $b$ \ dans \ $\mathbb{N}$ \ on a \ $F_{p,q}(a,b) =0 $ \ ainsi que \ $F_{p,\bar q}(a,b) = 0$ \ ce qui montre que la fonction consid\'er\'ee est \ $\mathscr{C}^{\infty}$ \ en \ $s = 0$. Ceci est \'evident \`a priori puisque l'on convole une fonction \ $\mathscr{C}^{\infty}$ \ avec une fonction localement int\'egrable \`a support compact.
       
       \bigskip
        
        \begin{cor}\label{Cas 1 complet}
        Dans la situation de la proposition pr\'ec\'edente, soient \ $j$ \ et \ $k$ deux entiers, et d\'efinissons pour \ $(a,b) \in V_{p,q}$ \ les fonctions
            \begin{align*}
       & F_{p,q}^{j,k}(a,b)[s] : = \frac{i}{4\pi}\int_{\vert u \vert \leq 1} \vert s-u\vert^{2a}.(s-u)^p.(Log\vert s-u\vert^2)^j.\vert u \vert^{2b}.u^q.(Log \vert u\vert^2)^k.du \wedge d\bar u \\
       & F_{p,\bar q}^{j,k}(a,b)[s] : =  \frac{i}{4\pi}\int_{\vert u \vert \leq 1} \vert s-u\vert^{2a}.(s-u)^p.(Log\vert s-u\vert^2)^j.\vert u \vert^{2b}.\bar u^q.(Log \vert u\vert^2)^k.du \wedge d\bar u
       \end{align*}
       Alors on a
           \begin{align*}
       & F^{j,k}_{p,q}(a,b)[s] = P^{j,k}_{p,q}(a,b)[Log\vert s\vert^2].s^{p+q}.\vert s\vert^{2(a+b+1)} + s^{p+q}.\Phi_{p,q}(a,b,s) \\
       &  F_{p,\bar q}^{j,k}(a,b)[s] = P^{j,k}_{p,\bar q}(a,b)[Log\vert s\vert^2].s^p.\bar s^q.\vert s\vert^{2(a+b+1)} + \bar s^{q-p}.\Phi_{p,\bar q}(a,b,s) 
       \end{align*}
       o\`u \ $P^{j,k}_{p,q}$ \ et \ $P^{j,k}_{p,\bar q}$ \ sont des polyn\^omes de degr\'e \ $j+k$ \ dont les coefficients d\'ependent holomorphiquement de \ $(a,b)$ \ et dont les coefficients dominants sont donn\'es par \ $ F_{p,q}(a,b) $ \ et \ $ F_{p,\bar q}(a,b)$ \ respectivement, et o\`u les fonctions \ $\Phi_{p,q}^{j,k}$ \ et \ $\Phi_{p,\bar q}^{j,k}$ \ sont obtenues via l'op\'eration \ $\frac{\partial^{j+k}}{\partial^ja \partial^kb}$ \ sur les fonctions \ $\Phi_{p,q}$ \ et \ $\Phi_{p.\bar q}$ \ de la proposition pr\'ec\'edente.
       \end{cor}
       
       \parag{Preuve} Il suffit d'appliquer l'op\'erateur diff\'erentiel \ $\frac{\partial^{j+k}}{\partial^ja \partial^kb}$ \ dans l'assertion de la proposition pr\'ec\'edente. $\hfill \blacksquare$
       
       \bigskip
       
       \parag{Compl\'ement}  Dans le cas o\`u \ $a$ \ (ou bien \ $b$) \ est dans \ $\mathbb{N}$, avec \ $(a,b) \in U_{p,q}$,  le terme de degr\'e \ $k+j$, c'est \`a dire le coefficient de   \\ 
       $s^{p+q}.\vert s\vert^{2(a+b+1)}.\big[Log\vert s\vert^2\big]^{j+k}$ \ (resp. de\ $s^p.\bar s^q.\vert s\vert^{2(a+b+1)}.\big[Log\vert s\vert^2\big]^{j+k}$)  est nul. \\
       Pour avoir le terme singulier dominant non nul on doit donc calculer le coefficient de \ $s^{p+q}.\vert s\vert^{2(a+b+1)}.\big[Log\vert s\vert^2\big]^{j+k-1}$ \ (resp. \ $s^p.\bar s^q.\vert s\vert^{2(a+b+1)}.\big[Log\vert s\vert^2\big]^{j+k-1}$). Un calcul simple donne que ce coefficient vaut 
       \begin{align*}
       &   \Big[ \frac{\partial F_{p,q}}{\partial a} +     \frac{\partial F_{p,q}}{\partial b}\Big](a,b)         \quad {\rm (resp.} \qquad  \Big[ \frac{\partial F_{p,\bar q}}{\partial a} +     \frac{\partial F_{p,\bar q}}{\partial b}\Big](a,b)).                                           
       \end{align*}
       Comme on a \ $\Gamma(z-k) = \frac{(-1)^k}{k!}\frac{1}{z} + holomorphe(z)$ \ pour \ $k \in \mathbb{N}$ \ et \ $z$ \ voisin de \ $0$, on constate que pour obtenir le coefficient cherch\'e il suffit de remplacer (si par exemple c'est \ $a$ \ qui est dans \ $\mathbb{N}$) le facteur \ $1/\Gamma(-a)$ \ dans l'expression de \ $F_{p,q}$ \ par le nombre \ $(-1)^{a+1}.a!$ \ (resp. dans l'expression de \ $F_{p,\bar q}$ ) . Ceci montre que ces coefficients sont non nuls (car \ $a$ \ et \ $b$ \ ne peuvent \^etre simultan\'ement dans \ $\mathbb{N}$ \ puisque \ $a + b+ 1 \not\in \mathbb{N}$ \ par hypoth\`ese).\\
       Donc dans ces cas les polyn\^omes en \ $Log\vert s\vert^2 $ \   \ $P^{j,k}_{p,q}(a,b)$ \ et \ $P^{j,k}_{p,\bar q}(a,b)$ \ sont de degr\'e exactement \ $j+k-1$ \ . $\hfill \square$
       
       \subsection{Le second cas.}
       
        \bigskip
        
        \begin{prop}[Le second cas : \ $a,b$ \ non entiers et \ $a+b+1 \in \mathbb{N}$.]\label{Cas 2}
        On suppose maintenant que \ $(a,b) \in U_{p,q}$ \ mais que   \ $a+b+1 $ \ est un entier. Alors on a
         \begin{align*}
       & F_{p,q}(a,b)[s] = \big[\tilde{F}_{p,q}(a,b).Log\vert s\vert^2 + c_{p,q}(a,b)\big]. s^{p+q}.\vert s\vert^{2(a+b+1)} \ + \  s^{p+q}.\Psi_{p,q}(a,b,s) \\
       &  F_{p,\bar q}(a,b)[s] = \big[\tilde{F}_{p,\bar q}(a,b).Log\vert s\vert^2 + c_{p,\bar q}(a,b)\big].s^p.\bar s^q.\vert s\vert^{2(a+b+1)} + \bar s^{q-p}.\Psi_{p,\bar q}(a,b,s) 
       \end{align*}
       o\`u les coefficients \ $\tilde{F}_{p,q}(a,b)$ \ et \ $\tilde{F}_{p,\bar q}(a,b)$ \ sont donn\'es par les formules suivantes
       \begin{align*}
       & \tilde{F}_{p,q}(a,b) =  \frac{(-1)^{a+b}}{\Gamma(-a).\Gamma(-b)}.\frac{\Gamma(a+p+1).\Gamma(b+q+1)}{\Gamma(a+b+2).\Gamma(a+b+p+q+2)} \\
       & \tilde{F}_{p,\bar q}(a,b) =  \frac{(-1)^{a+b}}{\Gamma(-a).\Gamma(-b)}.\frac{\Gamma(a+p+1).\Gamma(b+q+1)}{\Gamma(a+b+q+2).\Gamma(a+b+p+2)}
        \end{align*}
      o\`u les  fonctions \ $\Psi_{p,q}$ \ et \ $\Psi_{p,\bar q}$ \  sont \ $\mathscr{C}^{\infty}$ \ sur \ $U_{p,q} \times D$, holomorphes sur \ $U_{p,q}$ \ pour \ $s \in D$ \ fix\'e, et les fonctions \ $c_{p,q}$ \ et \ $c_{p,\bar q}$ \ sont holomorphes sur \ $U_{p,q}$.
       \end{prop}
       
       \parag{Remarque} Quand on suppose de plus que \ $a$ \ et \ $b$ \ ne sont pas entiers, les nombres 
        $\tilde{F}_{p,q}(a,b) $ \ et \ $\tilde{F}_{p,\bar q}(a,b)$ \ ne sont pas nuls. $\hfill \square$
        
        \bigskip
       
       \parag{Preuve} Elle est analogue au cas de la proposition \ref{Cas 1} sauf qu'il faut prendre en compte l'apparition du logarithme puisque le fait que \ $a+b+1$ \ soit entier oblige \`a rencontrer dans la somme l'int\'egrale \ $\int_3^{1/\vert s\vert} \frac{d\rho}{\rho}$. Posons \ $ r_j  = q+j+ 2(a+b+1) $, et reprenons le calcul de l'int\'egrale \ $I_j$ \ (voir (A) dans la preuve de la proposition \ref{Cas 1}). Le terme en \ $s^{p+q}.\vert s\vert^{2(a+b+1)} .Log\vert s\vert$ \ aura pour coefficient \ $- \gamma_{a,q+j}^{r_j}$. On obtiendra ainsi  que
       \begin{align*}
       &  \tilde{F}_{p,q}(a,b) = - \sum_{j=0}^p \ (-1)^j.C_p^j.\gamma_{a,q+j}^{r_j}\\
       & \quad \quad  = (-1)^{a+b}\frac{\Gamma(a+1)}{\Gamma(-a).\Gamma(-b).\Gamma(a+b+2)}.\sum_{j=0}^p (-1)^j.C_p^j.\frac{\Gamma(b+q+j+1)}{\Gamma(a+b+q+j+2)} \\
       & \quad \quad  =  \frac{(-1)^{a+b}}{\Gamma(-a).\Gamma(-b)}\frac{\Gamma(a+p+1).\Gamma(b+q+1)}{\Gamma(a+b+2).\Gamma(a+b+p+q+2)}
       \end{align*}
       d'apr\`es le lemme \ref{Formule}.\\
       
       \smallskip
       
       Le calcule analogue pour le coefficient \ $ \tilde{F}_{p,\bar q}(a,b) $ \ donne,          \begin{align*}
       &  \tilde{F}_{p,\bar q}(a,b) = - \sum_{j=0}^p \ (-1)^j.C_p^j.\gamma_{a,j-q}^{r_j}\\
       & \quad \quad =   (-1)^{a+b+q}\frac{\Gamma(a+1)}{\Gamma(-a).\Gamma(-b-q).\Gamma(a+b+q+2)}.\sum_{j=0}^p (-1)^j.C_p^j.\frac{\Gamma(b+j+1)}{\Gamma(a+b+j+2)} \\
          & \quad \quad  =  \frac{(-1)^{a+b+q}}{\Gamma(-a).\Gamma(-b-q)}\frac{\Gamma(a+p+1).\Gamma(b+1)}{\Gamma(a+b+q+2).\Gamma(a+b+p+2)}\\
          & \quad \quad = \frac{(-1)^{a+b}}{\Gamma(-a).\Gamma(-b)}\frac{\Gamma(a+p+1).\Gamma(b+q+1)}{\Gamma(a+b+q+2).\Gamma(a+b+p+2)}
       \end{align*}
       d'apr\`es la formule des compl\'ements. $\hfill \blacksquare$

       \bigskip
       
          \begin{cor}\label{Cas 2 complet}
        Dans la situation de la proposition pr\'ec\'edente, soient \ $j$ \ et \ $k$ deux entiers, et d\'efinissons pour \ $(a,b) \in U_{p,q}$ \ v\'erifiant \ $a+b+1 \in \mathbb{N}$,  les fonctions
            \begin{align*}
       & F_{p,q}^{j,k}(a,b)[s] : = \frac{i}{4\pi}\int_{\vert u \vert \leq 1} \vert s-u\vert^{2a}.(s-u)^p.(Log\vert s-u\vert^2)^j.\vert u \vert^{2b}.u^q.(Log \vert u\vert^2)^k.du \wedge d\bar u \\
       & F_{p,\bar q}^{j,k}(a,b)[s] : =  \frac{i}{4\pi}\int_{\vert u \vert \leq 1} \vert s-u\vert^{2a}.(s-u)^p.(Log\vert s-u\vert^2)^j.\vert u \vert^{2b}.\bar u^q.(Log \vert u\vert^2)^k.du \wedge d\bar u
       \end{align*}
       Alors on a
           \begin{align*}
       & F^{j,k}_{p,q}(a,b)[s] = P^{j,k}_{p,q}(a,b)[Log\vert s\vert^2].s^{p+q}.\vert s\vert^{2(a+b+1)} + s^{p+q}.\Phi_{p,q}(a,b,s) \\
       &  F_{p,\bar q}^{j,k}(a,b)[s] = P^{j,k}_{p,\bar q}(a,b)[Log\vert s\vert^2].s^p.\bar s^q.\vert s\vert^{2(a+b+1)} + \bar s^{q-p}.\Phi_{p,\bar q}(a,b,s) 
       \end{align*}
       o\`u \ $P^{j,k}_{p,q}$ \ et \ $P^{j,k}_{p,\bar q}$ \ sont des polyn\^omes de degr\'e \ $j+k+1$ \ dont les coefficients d\'ependent holomorphiquement de \ $(a,b)$ \ et dont les coefficients dominants sont donn\'es par \ $ \tilde{F}_{p,q}(a,b) $ \ et \ $ \tilde{F}_{p,\bar q}(a,b)$ \ respectivement, et o\`u les fonctions \ $\Phi_{p,q}^{j,k}$ \ et \ $\Phi_{p,\bar q}^{j,k}$ \ sont obtenues via l'op\'eration \ $\frac{\partial^{j+k}}{\partial^ja \partial^kb}$ \ sur les fonctions \ $\Phi_{p,q}$ \ et \ $\Phi_{p.\bar q}$ \ de la proposition pr\'ec\'edente.
       \end{cor}
       
         \parag{Preuve} Il suffit \`a nouveau  d'appliquer l'op\'erateur diff\'erentiel \ $\frac{\partial^{j+k}}{\partial^ja \partial^kb}$ \ dans l'assertion de la proposition pr\'ec\'edente. $\hfill \blacksquare$

\bigskip

Le dernier cas \`a traiter est celui o\`u \ $a$ \ et \ $b$ \ sont entiers. Ceci ne peut s'obtenir comme pr\'ec\'edemment par d\'erivation du cas o\`u les logarithmes n'apparaissent pas. Il faut donc traiter directement l'analogue des corollaires \ref{Cas 1 complet} et \ref{Cas 2 complet}.

\subsection{Le dernier cas.}

\bigskip

       \begin{prop}[Le cas \ $a$ \ et \ $b$ \ entiers.] \label{Cas 3 complet} Donnons-nous deux entiers \\
        $0 \leq p \leq q$. Supposons  maintenant que \ $a$ \ et \ $b$ \ sont dans \ $\mathbb{Z}$ \ et que l'on a \ $a+p/2 > -1$ \ et \ $b+ q/2 > - 1$. Pour \ $(j,k) \in (\mathbb{N}^*)^2$ \ et \ $s \in D^*$ \  posons alors
       \begin{align*}
      & F^{j,k}_{p,q}(a,b)[s] : = \frac{i}{4\pi} \int_{\vert u \vert \leq 1} \vert s-u\vert^{2a}.(s-u)^p.(Log\vert s-u\vert^2)^j.\vert u\vert^{2b}.u^q.(Log\vert u\vert^2)^k.du\wedge d\bar u \\
        & F^{j,k}_{p,\bar q}(a,b)[s] : = \frac{i}{4\pi} \int_{\vert u \vert \leq 1} \vert s-u\vert^{2a}.(s-u)^p.(Log\vert s-u\vert^2)^j.\vert u\vert^{2b}.\bar u^q.(Log\vert u\vert^2)^k.du\wedge d\bar u 
        \end{align*}
          Alors on a
           \begin{align*}
       & F^{j,k}_{p,q}(a,b)[s] = \check{P}^{j,k}_{p,q}(a,b)[Log\vert s\vert^2].s^{p+q}.\vert s\vert^{2(a+b+1)} + s^{p+q}.\Phi_{p,q}(a,b)[\vert s\vert^2] \\
       &  F^{j,k}_{p,\bar q}(a,b)[s] = \check{P}^{j,k}_{p,\bar q}(a,b)[Log\vert s\vert^2].s^p.\bar s^q.\vert s\vert^{2(a+b+1)} + \bar s^{q-p}.\Phi_{p,\bar q}(a,b)[\vert s\vert^2] \tag{C}
       \end{align*}
       o\`u \ $\check{P}^{j,k}_{p,q}$ \ et \ $\check{P}^{j,k}_{p,\bar q}$ \ sont des polyn\^omes de degr\'e exactement \ $j+k-1$ \   et o\`u les fonctions \ $\Phi_{p,q}^{j,k}(a,b)$ \ et \ $\Phi_{p,\bar q}^{j,k}(a,b)$ \ sont analytiques r\'eelles.
        \end{prop}
        
        \parag{Preuve} Commen{\c c}ons par rappeler que pour \ $j =0$ \  ou \ $k = 0$ \ (cas exclus de l'\'enonc\'e ci-dessus) les fonctions consid\'er\'ees sont \ $\mathscr{C}^{\infty}$ \ comme convol\'ees d'une fonction \ $\mathscr{C}^{\infty}$ \ et d'une fonction localement int\'egrable \`a support compact. \\
        Le changement de variable \ $ u = t.s$, pour \ $s \not= 0$ \ fix\'e donne
        \begin{align*}
        & F^{j,k}_{p,q}(a,b)[s] = s^{p+q}.\vert s\vert^{2(a+b+1)}.I(s)
        \end{align*}
       o\`u nous avons pos\'e
       \begin{align*}   
        & I(s) : =  \int_{\vert t \vert \leq 1/\vert s\vert} \vert 1-t\vert^{2a}.(1-t)^p.(Log\vert 1-t\vert^2 + Log\vert s\vert^2)^j.\vert t\vert^{2b}.t^q.(Log\vert s.t\vert^2)^k.dt\wedge d\bar t 
        \end{align*}
        ainsi qu'une expression analogue pour \ $ F^{j,k}_{p,\bar q}(a,b)[s] $. Des calculs analogues \`a ceux d\'ej\`a d\'etaill\'es plus haut montrent facilement que l'on a des expressions du type \ $({\rm C})$, mais avec des polyn\^omes en \ $Log\vert s\vert^2$ \ \`a priori de degr\'es  inf\'erieurs ou \'egaux \`a \ $k+j+1$.\\
        Nous allons d\'emontrer l'assertion sur le degr\'e de ces polyn\^omes par r\'ecurrence sur \ $j+k = n \geq 2$. 
        Pour \ $n = 2$ \ on a n\'ecessairement \ $j = k = 1$ \ et il s'agit de montrer que les polyn\^omes \ $\check{P}^{1,1}_{p,q}(a,b)$ \ et \ $\check{P}^{1,1}_{p,q}(a,b)$ \ sont de degr\'es exactement  \'egal \`a \ $1$. Le pas de r\'ecurrence qui va suivre montrera qu'ils sont de degr\'es au plus \'egal \`a \ $1$. Il nous suffit donc de montrer que le coefficient de \ $Log\vert s\vert^2$ \ est non nul.\\
       Ceci r\'esulte du calcul de la constante \ $\gamma_{1,1}$ \ qui est fait au paragraphe 3.5.
               
        \bigskip

        Supposons  d\'emontr\'e que pour \ $n \geq 2$ \  le degr\'e des polyn\^omes \ $\check{P}^{j,k}_{p,q}(a,b)$ \ et \ $\check{P}^{j,k}_{p,\bar q}(a,b)$ \ est au plus \'egal \`a \ $j+k-1 $ \ et montrons ceci pour un couple \ $(j,k) \in (\mathbb{N}^*)^2$ \ v\'erifiant \ $j+k = n+1$. Soit \ $\lambda \in \mathbb{C}^*$ \ et calculons la diff\'erence
        $$ F^{j,k}_{p,q}(a,b)[\lambda.s] - \lambda^{p+q}.\vert \lambda\vert^{2(a+b+1)}.F^{j,k}_{p,q}(a,b)[s] $$
        en utilisant le changement de variable \ $u = \lambda.v$. On obtient \ $  \frac{i}{4\pi}.\lambda^{p+q}.\vert \lambda\vert^{2(a+b+1)}$ \ multipli\'e par
        \begin{align*}
        &\int_{1\leq \vert v\vert \leq 1/\vert \lambda\vert} \vert s - v\vert^{2a}.(s-v)^p.(Log\vert \lambda\vert + Log\vert s-v\vert)^j.\vert v\vert^{2b}.v^q.(Log\vert \lambda\vert + Log\vert v\vert)^k.dv\wedge d\bar v  \\
        &\qquad \qquad \qquad \qquad + \int_{\vert v\vert \leq 1} \vert s-v\vert^{2a}.(s-v)^p.\vert v\vert^{2b}.v^q.Z.dv\wedge d\bar v 
\\ 
        &{\rm avec} \quad Z : = \Big[(Log\vert \lambda\vert + Log\vert s-v\vert)^j.(Log\vert \lambda\vert + Log\vert v\vert)^k - ( Log\vert s-v\vert)^j.( Log\vert v\vert)^k\Big]    
            \end{align*}
            On constate que la premi\`ere int\'egrale est \ $\mathscr{C}^{\infty}$ \ sur \ $D$, et que la seconde est une combinaison lin\'eaire des fonctions \ $F^{j',k'}_{p,q}(a,b)$ \ avec \ $j'+k' \leq n $. On en d\'eduit facilement notre assertion.\\
            De plus, si \ $\gamma(j,k)$ \ d\'esigne le coefficient de \ $(Log\vert s\vert)^{j+k-1}$ \ dans \ $ \check{P}^{j,k}_{p,q}(a,b)$ \ le calcul ci-dessus donne facilement la relation 
            $$\gamma(j,k) = (j+k)! \gamma(1,1) \quad \forall j, k \geq 1$$
      ce qui montre que si \ $\gamma(1,1)$ \ est non nul, il en est de m\^eme pour tous les \\
       $\gamma(j,k), \forall (j,k) \in (\mathbb{N}^*)^2 $. Nous allons montrer  au paragraphe suivant que la constante  \ $\gamma(1,1) $ \ (qui d\'epend de \ $a,b,p,q$)  est non nulle, ce qui  ach\`evera la preuve. $\hfill \blacksquare$
            
            \bigskip
                        
   \subsection{Le calcul de \ $\gamma(1,1)$.}

       \bigskip

       \begin{lemma}\label{Calcul 1}
       Pour \ $0 \leq x < 1 $ \ et \ $p \in \mathbb{Z}$ \ on a
       \begin{align*}
       & C_p(x) : =  \frac{1}{2\pi}\int_0^{2\pi} Log\vert 1- x.e^{i\theta}\vert^2.e^{ip\theta}.d\theta = \frac{x^{\vert p\vert}}{\vert p\vert} \quad {\rm pour} \quad p \not= 0 \\
       & C_0(x) :  =  \frac{1}{2\pi}\int_0^{2\pi} Log\vert 1- x.e^{i\theta}\vert^2.d\theta = 0.
       \end{align*}
       Pour \ $x > 1$ \ on a \ $C_p(x) = C_p(1/x)$ \ pour \ $p \not= 0$ \ et \ $C_0(x) = Log\, x^2 $.
       \end{lemma}
       
       \parag{Preuve} Posons pour \ $x \in ]0,1[$ 
        $$A_p : = \frac{1}{2i\pi}\int_{\vert z\vert = x} Log(1-z).z^p \frac{dz}{z} .$$
        Alors on a \ $C_p(x) = x^{-p}.A_p - x^p.\bar A_{-p}$. Comme la formule de Cauchy donne
        \begin{align*}
        & A_p = 0 \quad {\rm pour} \quad  p \geq 0  \quad {\rm et} \\
        & A_p = -\frac{1}{p} \quad {\rm pour} \quad p < 0
        \end{align*}
        on conclut facilement.$\hfill \blacksquare$
        
        \begin{prop}\label{Calcul 2}
        Pour \ $p,q$ \ deux entiers naturels, posons pour \ $s \in D$
        $$ F_{p,q}(s) : = \frac{i}{4\pi}\int_{\vert u \vert \leq 1} \ u^p.\bar u^q.Log\vert s-u\vert.Log\vert u\vert.du\wedge d\bar u .$$
        Alors le coefficient du terme en \ $s^p.\bar s^q.Log\vert s\vert^2 $ \ dans le d\'eveloppement asymptotique en \ $s = 0$ \ de \ $F_{p,q}$ \ vaut \ $- \frac{1}{4(p+1)(q+1)}$.
                \end{prop}
                
                On remarquera que le terme en  \ $s^p.\bar s^q.Log\vert s\vert^2 $ \  est le seul terme non \ $\mathscr{C}^{\infty}$ \ dans le d\'eveloppement de cette fonction \`a l'origine.
         
         \parag{Preuve} Commen{\c c}ons par le cas \ $p \not= q$. On a, en posant \ $u = s.t$ \ pour \ $s \not= 0$
         \begin{align*}
         & F_{p,q}(s) = s^{p+1}.\bar s^{q+1}\frac{i}{4\pi}\int_{\vert t\vert \leq 1/\vert s\vert} t^p.\bar t^q.Z.dt\wedge d\bar t  \\
         & {\rm avec} \quad  Z : = (Log\vert s\vert + Log\vert 1-t\vert).(Log\vert s\vert + Log\vert t\vert)
         \end{align*}
         Cela donne les trois termes suivants
         \begin{align*}
         & A : =  s^{p+1}.\bar s^{q+1}.(Log\vert s\vert)^2.I_1 \\
         & B : =   s^{p+1}.\bar s^{q+1}.(Log\vert s\vert).I_2 \\
         & C : =  s^{p+1}.\bar s^{q+1}.I_3
         \end{align*}
         o\`u les int\'egrales \ $I_1, I_2$ \ et \ $I_3$ \ vont \^etre examin\'ees ci-dessous.\\
         Remarquons d\'ej\`a que le d\'eveloppement asymptotique de \ $A$ \ ne donnera jamais de contribution au terme qui nous int\'eresse.\\
         Pour  \ $B$ \ nous cherchons le terme constant dans le d\'eveloppement asymptotique de 
         $$ I_2(s) : = \frac{i}{4\pi}\int_{\vert t\vert \leq 1/\vert s\vert} t^p.\bar t^q.(Log\vert 1-t\vert +  Log\vert t\vert).dt\wedge d\bar t  \quad.$$
         Cherchons d\'ej\`a le terme constant dans le d\'eveloppement de l'int\'egrale
         \begin{align*}
      &   I'_2(s) : =  \frac{i}{4\pi}\int_{1 \leq \vert t\vert \leq 1/\vert s\vert} t^p.\bar t^q.(Log\vert 1-t\vert +  Log\vert t\vert).dt\wedge d\bar t \\
                 & \quad = \frac{1}{2}.\int_1^{1/\vert s\vert} \rho^{p+q+1}.C_{p-q}(\rho).d\rho \ + \    \int_1^{1/\vert s\vert} \rho^{p+q+1} Log\, \rho .d\rho \quad .
         \end{align*}
         Comme le lemme pr\'ec\'edent donne \ $C_{p-q}(\rho) = \frac{\rho^{-\vert p-q\vert}}{\vert p-q\vert}$ \ puisque l'on suppose \ $p \not= q$, on obtient facilement que le terme constant du d\'eveloppement de \ $I'_2(s)$ \ vaut 
         $$ -\frac{1}{2}.\frac{1}{p+q+2-\vert p-q\vert} \  + \ \frac{1}{(p+q+2)^2} .$$
         Il nous reste encore \`a \'evaluer la constante
         $$ I_2" : = \frac{i}{4\pi}\int_{\vert t\vert \leq 1} t^p.\bar t^q.(Log\vert 1-t\vert +  Log\vert t\vert).dt\wedge d\bar t $$
         ce qui est simple \`a l'aide du lemme pr\'ec\'edent : il donne
         \begin{align*}
         &  I_2"  =  \frac{1}{2}.\int_0^1 \rho^{p+q+1}.C_{p-q}(\rho).d\rho \ + \ \int_0^1  \rho^{p+q+1} Log\, \rho .d\rho \\
         & \quad =  \frac{1}{2.\vert p-q\vert}.\Big[\frac{\rho^{p+q+2+\vert p-q\vert}}{p+q+2+\vert p-q\vert}\Big]_0^1 \ - \  \frac{1}{(p+q+2)^2}\\
         & \quad = \frac{1}{2.\vert p-q\vert}.\frac{1}{p+q+2+\vert p-q\vert} - \frac{1}{(p+q+2)^2} 
         \end{align*}
      On trouve finalement, comme contribution de \ $I_2$ \ la constante
      \begin{align*}
      &  -\frac{1}{2\vert p-q\vert}.\Big[\frac{1}{p+q+2-\vert p-q\vert} \  -  \frac{1}{p+q+2+\vert p-q\vert} \Big]\\
            & \quad = -\frac{1}{4(p+1)(q+1)}
      \end{align*}
         
        Cherchons la contribution de \ ${\rm C}$ \ c'est \`a dire le terme en \ $Log \vert s\vert$ \ dans le d\'eveloppement asymptotique de 
       
       \begin{align*}
       & I_3(s) : = \frac{i}{4\pi}\int_{\vert t\vert \leq 1/\vert s\vert} t^p.\bar t^q.(Log\vert 1-t\vert).(Log\vert t\vert).dt\wedge d\bar t 
       \end{align*}
       Comme le terme constant ne nous int\'eresse pas, on peut se contenter de regarder
       \begin{align*}
      &  I'_3(s) : = \frac{i}{4\pi}\int_{1 \leq \vert t\vert \leq 1/\vert s\vert} t^p.\bar t^q.(Log\vert 1-t\vert).(Log\vert t\vert).dt\wedge d\bar t \\
      & \quad =  \int_1^{1/\vert s\vert} \rho^{p+q+1}.C_{p-q}(\rho).Log\,\rho.d\rho \\
       & \quad  =  \frac{1}{\vert p-q\vert}.\int_1^{1/\vert s\vert} \rho^{p+q+1-\vert p-q\vert}.Log\,\rho.d\rho 
      \end{align*}
          et il n'y a pas de terme en \ $Log \vert s\vert$ \ dans le d\'eveloppement asymptotique de \ $I_3(s)$, puisque \ $p+q+2-\vert p-q\vert \geq 2$.\\
          Le cas \ $p = q$ \ est analogue en utilisant le calcul de \ $C_0(x)$ \ dans le lemme pr\'ec\'edent, en prenant garde au cas \ $x > 1$. $\hfill \blacksquare$
          
          \begin{cor}\label{Calcul 3}
          Soient \ $a,b,p,q$ \ des entiers naturels. Pour \ $s \in D^*$, posons
          \begin{align*}
          &  F_{p,q}(a,b)[s] : = \frac{1}{4i\pi}\int_{\vert u \vert \leq 1} \vert s-u\vert^{2a}.(s-u)^p.Log\vert s-u\vert.\vert u\vert^{2b}.u^q.Log\vert u\vert.du\wedge d\bar u \\
           &  F_{p,\bar q}(a,b)[s] : = \frac{1}{4i\pi}\int_{\vert u \vert \leq 1} \vert s-u\vert^{2a}.(s-u)^p.Log\vert s-u\vert.\vert u\vert^{2b}.\bar u^q.Log\vert u\vert.du\wedge d\bar u
           \end{align*}
           Le coefficient du terme en \ $\vert s\vert^{2(a+b+1)}.s^{p+q}.Log\vert s\vert $ \ dans le d\'eveloppement asymptotique de \ $ F_{p,q}(a,b)$ \ en \ $s = 0$ \ est \'egal \`a 
           \begin{align*}
           & -\frac{1}{4}\frac{\Gamma(a+1).\Gamma(b+1)}{\Gamma(a+b+2)}.\frac{\Gamma(a+p+1).\Gamma(b+q+1)}{\Gamma(a+b+p+q+2)} \quad .
           \end{align*}
            Le coefficient du terme en \ $\vert s\vert^{2(a+b+1)}.s^p.\bar s^q.Log\vert s\vert $ \ dans le d\'eveloppement asymptotique de \ $ F_{p,\bar q}(a,b)$ \ en \ $s = 0$ \ est \'egal \`a 
            $$ -\frac{1}{4}\frac{\Gamma(a+1).\Gamma(b+q+1)}{\Gamma(a+b+q+2)}.\frac{\Gamma(a+p+1).\Gamma(b+1)}{\Gamma(a+b+p+2)} \quad .$$
            \end{cor}
            
            On remarquera \`a nouveau que, pour chacune de ces fonctions,  le terme consid\'er\'e  est le seul terme non \ $\mathscr{C}^{\infty}$ \ du d\'eveloppement asymptotique.
            
            \parag{Preuve} La formule du bin\^ome et la proposition pr\'ec\'edente donne que le coefficient cherch\'e vaut, pour la fonction \ $ F_{p,q}(a,b)$,
            \begin{align*}
            & \frac{-1}{4}. \sum_{j=0}^{a+p}\sum_{k=0}^a (-1)^{j+k}.C_{a+p}^j.C_a^k \frac{1}{(b+q+j+1)(b+k+1)  } \quad {resp.} \\
            &  \frac{-1}{4}. \sum_{j=0}^{a+p}\sum_{k=0}^a (-1)^{j+k}.C_{a+p}^j.C_a^k\frac{1}{(b+j+1)(b+q+k+1) }
            \end{align*}
            En utilisant la formule
            $$ \sum_{k=0}^m (-1)^k.C_m^k\frac{1}{n+k} = \frac{\Gamma(m+1).\Gamma(n)}{\Gamma(m+n+1)}$$
            on obtient facilement le r\'esultat annonc\'e. L'autre cas est analogue. $\hfill \blacksquare$
            
            \parag{Remarque} Le corollaire pr\'ec\'edent montre, en particulier,  que ces coefficients ne sont jamais nuls.

     \newpage

 \section{R\'ef\'erences.}
    
    \begin{itemize}
      \item{[B.\,82]} Barlet, D. \textit{D\'eveloppements asymptotiques des fonctions obtenues par int\'egration sur les fibres},  Inv. Math. vol. 68 (1982), p. 129-174.
        \item{[B.\,86]}  Barlet, D. \textit{Calcul de la forme hermitienne canonique pour \ $X^a + Y^b + Z^c $}, in Sem. P. Lelong, Lecture Notes, vol. 1198 Springer Verlag (1986), p. 35-46.
        \item{[Bj.\,93]} Bjork,J.-E. {\it Analytic D-Modules and Applications}, Kluwer Academic Publishers, Dordrecht/Boston/London 1993.
        \item{[B.-S.\,74]} Brian{\c c}on, J. et Skoda, H. {\it Sur la cl\^oture int\'egrale d'un id\'eal de germes de fonctions holomorphes en un point de \ $\mathbb{C}^n$}, C.R.Acad.Sci. Paris s\'erie A, 278 (1974), p.949-951.
        \item{[K.\,76]} Kashiwara,M. {\it b-function and holonomic systems}, Inv. Math. 38, (1976), p.33-53.
        \item{[M.\,74]} Malgrange, B.  {\it Int\'egrale asymptotique et monodromie.} Ann. Scient. Ec. Norm. Sup. , t.7 (1974), p.405-430.
        \item{[M.\,83]} Malgrange, B. {\it Polyn\^ome de Bernstein-Sato et cohomologie \'evanescente}, Ast\'erisque 101-102 (1983), p.243-267.
        \item{[Sak.\,73]} Sakamoto, K. {\it Milnor fibering and their characteristic maps}, Proc. Intern. Conf. on manifolds and Related Topics in Topology, Tokyo 1973.
        \item{[S.-T.\,71]} Sebastiani, M. and Thom, R. {\it Un r\'esultat sur la monodromie}, Inv. Math. 13 (1971), p. 90-96.
        \end{itemize}
        
        \bigskip

Barlet Daniel, Institut Elie Cartan UMR 7502  \\
Nancy-Universit\'e, CNRS, INRIA  et  Institut Universitaire de France, \\
BP 239 - F - 54506 Vandoeuvre-l\`es-Nancy Cedex.France.\\
e-mail : barlet@iecn.u-nancy.fr.

\end{document}